\documentclass[a4paper,12pt]{amsart}
\usepackage{amsmath, amssymb, amsthm}
\usepackage{amssymb}

\newcommand{\K}{\mathcal{K}}

\newcommand{\rank}{\mathrm{Rank}\,}
\newcommand{\chii}{\raisebox{2pt}{$\chi$}}
\newcommand{\C}{\mathbb{C}}
\newcommand{\R}{\mathbb{R}}
\newcommand{\tens}{\otimes}
\newcommand{\dsum}{\oplus}

\newcommand{\fixlist}{\mbox{} \@afterheading}

\newcommand{\iso}{\cong}
\newcommand{\dunion}{\amalg}

\renewcommand{\setminus}{\backslash}
\renewcommand{\emptyset}{\varnothing}

\newcommand{\labelledthing}[2]{\hspace{4pt}\buildrel {#2} \over #1 \hspace{3pt}} 
\newcommand{\labelledleftarrow}{\labelledthing{\longleftarrow}}
\newcommand{\labelledrightarrow}{\labelledthing{\longrightarrow}}
\usepackage{hyperref}

\newcounter{theorem}
\numberwithin{theorem}{section}
\newtheorem{thm}[theorem]{Theorem}
\newtheorem{lem}[theorem]{Lemma}
\newtheorem{prop}[theorem]{Proposition}
\newtheorem{cor}[theorem]{Corollary}
\newtheorem{defn}[theorem]{Definition}
\newtheorem{rmk}[theorem]{Remark}

\theoremstyle{remark}
\newtheorem*{remark*}{Remark}

\numberwithin{equation}{section}

\newcommand{\alabel}{\label}

\newcommand{\Cu}{\mathcal{C}u}
\newcommand{\culeq}{\preceq_{Cu}}
\newcommand{\cueq}{\sim_{Cu}}
\newcommand{\mvnleq}{\preceq}
\newcommand{\mvneq}{\sim}
\newcommand{\hered}[1]{\mathrm{Her}\left(#1\right)}

\newcommand{\e}{\epsilon}
\newcommand{\dl}{\delta}
\newcommand{\jsZ}{\mathcal{Z}}
\newcommand{\id}{\mathrm{id}}
\newcommand{\her}[1]{\mathrm{Her}\left(#1\right)}
\newcommand{\ev}{\mathrm{ev}}
\newcommand{\ccite}[2]{\cite[#1]{#2}}
\newcommand{\sqbox}{\qedhere}

\begin{document}

\title[The Cuntz semigroup of $C_0(X,A)$]{The Cuntz semigroup of continuous functions into certain simple $C^*$-algebras}
\author{Aaron Tikuisis}
\thanks{The author was supported by an NSERC CGS-D scholarship}


\maketitle

\begin{abstract}
This paper contains computations of the Cuntz semigroup of separable $C^*$-algebras of the form $C_0(X,A)$, where $A$ is a unital, simple, $\jsZ$-stable ASH algebra.
The computations describe the Cuntz semigroup in terms of Murray-von Neumann semigroups of $C(K,A)$ for compact subsets $K$ of $X$.
In particular, the computation shows that the Elliott invariant is functorially equivalent to the invariant given by the Cuntz semigroup of $C(\mathbb{T},A)$.
These results are a contribution towards the goal of using the Cuntz semigroup in the classification of well-behaved non-simple $C^*$-algebras.

\end{abstract}

\section{Introduction\alabel{intro}}

The Cuntz semigroup is an isomorphism invariant for $C^*$-algebras, consisting of an ordered semigroup that is constructed using positive elements of a $C^*$-algebra in much the same way that the Murray-von Neumann semigroup is constructed using projections.
Whereas Murray-von Neumann comparison theory tells us a lot about the structure of von Neumann algebras, $C^*$-algebras generally have much fewer projections than von Neumann algebras, and so by using positive elements, the Cuntz semigroup has the capability to detect significantly more structure.
This sensitivity of the Cuntz semigroup makes it an excellent candidate to distinguish $C^*$-algebras, both simple and non-simple.
For simple, well-behaved algebras, it contains the same information as some classical invariants \cite{BrownPereraToms}, which are a major part of the invariant used in the classification of large classes of $C^*$-algebras \cite{ElliottGongLi:AHclassification,Lin:AsympClassification,LinNiu:Lifting}.
In the non-simple case, the classical invariants are insufficient, whereas the Cuntz semigroup naturally handles ideal structure; it has already been used for non-simple classification in some promising first steps \cite{CiupercaElliottSantiago,Robert:NCCW}.
Yet the sensitivity of the Cuntz semigroup also makes it difficult to compute or describe in any concrete terms, except with the simple well-behaved algebras upon which it is equivalent to a better-studied invariant.

The main result of this paper consists of a computation of the Cuntz semigroup of separable algebras $C_0(X,A)$ when $A$ is a simple, unital, $\jsZ$-stable approximately subhomogeneous (ASH) algebra.
The class of ASH algebras includes many important, naturally occurring examples from dynamical systems \cite{LinPhillips:LimitDecomp,Phillips:Tori,Thomsen:Solenoids,TomsWinter:MinClassification}, as well as all approximately finite (AF) $C^*$-algebras and the Jiang-Su algebra $\jsZ$.
In fact, there is no known example of a simple, separable, stably finite, nuclear $C^*$-algebra which is not ASH.
It should be noted that, for a simple, unital ASH algebra, the regularity property of being $\jsZ$-stable is equivalent to having slow dimension growth and being non-type I (see Proposition \ref{ASH-SlowCharacterization}).

A key component of the Cuntz semigroup computation is the observation that, for a positive element $a$ of such $C_0(X,A)$, its Cuntz class is determined by its value $d_\tau(a)$ under every dimension function $d_\tau$ on $C_0(X,A)$, along with the Murray-von Neumann class of any image of $a$ in any quotient for which the image of $a$ is equivalent to a projection (such a quotient is necessarily $C(K,A)$ where $K \subseteq X$ is compact).

The Murray-von Neumann semigroups of $C(K,A)$ enter the description of the Cuntz semigroup of the algebras $C_0(X,A)$ via the object
\[ V_c(Y,A) := \{\text{projection valued functions } p \in C_b(Y, A \tens \K)\}/ \sim_c \]
where $p \sim_c q$ if for every compact set $K$, there exists $v \in C(K,A \tens \K)$ such that $p|_K = vv^*$ and $v^*v = q|_K$ (ie.\ $p|_K$ is Murray-von Neumann equivalent to $q|_K$), and where $Y$ is the difference of two open subsets of $X$.
We use $\langle p \rangle$ to denote the equivalence class in $V_c(Y,A)$ of the projection-valued function $p:Y \to A \tens \K$.
Of course, if $Y$ is compact then $V_c(Y,A)$ is the same as the Murray-von Neumann semigroup $V(C(Y,A))$, and in general, it is the inverse limit of $V(C(K,A))$ over compact subsets $K$ of $Y$ (Corollary \ref{VcDescr}).
The computation result may now be stated:

\begin{thm}\alabel{MainThm}
Let $A$ be a simple, unital $\jsZ$-stable ASH algebra and let $X$ be a second countable, locally compact Hausdorff space.
Then $\Cu(C_0(X,A))$ may be identified with pairs $\left(f,\left(\langle q_{[p]} \rangle\right)_{[p] \in V(A)}\right)$, where 
\begin{itemize}
\item $f:X \to \Cu(A)$ is a function which is lower semicontinuous with respect to $\ll$, and 
\item for each $[p] \in V(A)$, $\langle q_{[p]} \rangle$ is an element of $V_c(f^{-1}([p]),A)$ such that $[q_{[p]}(x)] = [p]$ in $V(A)$ for each $x \in f^{-1}([p])$.
\end{itemize}
The ordering is given by $\left(f,\left(\langle q_{[p]} \rangle\right)_{[p] \in V(A)}\right) \leq \left(g, \left(\langle r_{[p]} \rangle\right)_{[p] \in V(A)}\right)$ if $f(x) \leq g(x)$ for each $x$, and for each $[p] \in V(A)$,
\[ \langle q_{[p]}|_{f^{-1}([p]) \cap g^{-1}([p])} \rangle = \langle r_{[p]}|_{f^{-1}([p]) \cap g^{-1}([p])} \rangle. \]
The addition is given by 
\[ \left(f,\left(\langle q_{[p]} \rangle\right)_{[p] \in V(A)}\right) + \left(g, \left(\langle r_{[p]} \rangle\right)_{[p] \in V(A)}\right) = \left(f+g, \left(\langle s_{[p]} \rangle\right)_{[p] \in V(A)}\right), \]
where for every pair of projections $0 \leq p' \leq p \in A \tens \K$, we have
\[ s_{[p]}|_{f^{-1}([p']) \cap g^{-1}([p-p'])} = q_{[p']} + r_{[p-p']}. \]
(We have that $(f+g)^{-1}([p])$ breaks into disjoint components $f^{-1}([p']) \cap g^{-1}([p-p'])$, and so this definition of $s_{[p]}$ is continuous.)
\end{thm}

This result provides a start to the problem of finding the Cuntz semigroup for general well-behaved non-simple $C^*$-algebras.
(The adjective ``well-behaved'' is intentionally vague, referring to a variety of regularity conditions, (A)-($\Delta$) in the introduction of \cite{Winter:drZstable}; the algebras studied here have the regularity property of being $\jsZ$-stable.)
Although the description in Theorem \ref{MainThm} is not short, it is as nice as it may be, in the sense that for any $C^*$-algebra $A$, Cuntz comparison of positive elements in $C_0(X,A)$ always takes into account the data that appears in this computation.
Sometimes, the data appearing here is not even enough to determine the Cuntz class, as shown by Leonel Robert and the author in \ccite{Example 2, in 6.1}{CommutativeCuntz} for $X=S^4$ and $A=\C$.
Theorem \ref{MainThm} attests to the notion that the algebras $C_0(X,A)$ appearing in the theorem are well-behaved, yet also indicates how complex the Cuntz semigroup of well-behaved, non-simple $C^*$-algebras may be (as compared to the case with simple, well-behaved $C^*$-algebras).

Strikingly, the regularity of the algebra $A$ wins over the potential topological irregularity of the space $X$, since the dimension of $X$ is entirely unrestricted.
Leonel Robert and the author proved in \ccite{Section 5.1}{CommutativeCuntz} that the description in Theorem \ref{MainThm} also holds for $A=\C$ when the space $X$ is restricted to have covering dimension at most three.
Many of the arguments here build on those used in \cite{CommutativeCuntz}.

As an application of Theorem \ref{MainThm}, by using $X = \mathbb{T}$, we show in Section \ref{Ell-CuT-Eq} that the Elliott invariant is equivalent to the invariant $\Cu_\mathbb{T}(A) := (\Cu(C(\mathbb{T},A)), [1]_{\Cu(C(\mathbb{T},A))})$.
It follows that the Elliott invariant and $\Cu_\mathbb{T}(\cdot)$ classify exactly the same classes of simple, unital, $\jsZ$-stable ASH algebras.

In Section \ref{Prelim}, we recall basic facts and definitions concerning Cuntz and Murray-von Neumann comparison, along with the semigroups that are built from each of these comparison relations.
Section \ref{GeneralConsid} contains general considerations about the Cuntz semigroup of $C_0(X,A)$ where $A$ is a stably finite $C^*$-algebra; specifically, we study the map $x \mapsto [a(x)]$ from $X$ to $\Cu(A)$ (for a fixed element $a \in C_0(X, A \tens \K)_+$) and the semigroup $V_c(X,A)$.
In Section \ref{ComparisonSection}, we prove half of the main theorem; namely, we show that a Cuntz equivalence invariant (giving rise to the data appearing in Theorem \ref{MainThm}) is complete, ie.\ it completely determines the Cuntz order.
The other half of the proof of the main theorem consists of describing the range of this invariant.
Here, we show that every element of the semigroup described in Theorem \ref{MainThm} occurs as the value of the invariant on some positive element; this part is shown in Section \ref{DataAttained}, and a preliminary result there, Proposition \ref{Cu-Interpolation}, states that the Cuntz semigroup of a simple, $\jsZ$-stable algebra has Riesz interpolation.
In Section \ref{Ell-CuT-Eq}, we describe the functorial equivalence of the Elliott invariant and the Cuntz semigroup of $C(\mathbb{T},A)$.
Section \ref{OldCu} contains an analogue of Theorem \ref{MainThm} for the original Cuntz semigroup $W(C_0(X,A))$, and Section \ref{MvN} contains some further remarks about the Murray-von Neumann semigroups of $C(K,A)$ for $A$ as in Theorem \ref{MainThm}.
\section{Preliminaries\alabel{Prelim}}

\subsection{The semigroups of Cuntz and Murray-von Neumann}
The following pre-order relation was defined in \cite{Cuntz:dimensionfunctions}, and is the first step in constructing the Cuntz semigroup.

\begin{defn}
Let $A$ be a $C^*$-algebra.
The pre-order $\culeq$ is defined on $A_+$ by the following.
For positive elements $a,b \in A$, we have $a \culeq b$ if there exists a sequence $(x_n) \in A$ such that
\[ a = \lim_{n \to \infty} x_n^*bx_n. \]
The order $\cueq$ is given by $a \cueq b$ if $a \culeq b$ and $b \culeq a$.
\end{defn}

The Cuntz semigroup, as defined here and denoted $\Cu(A)$, is different from the original definition in \cite{Cuntz:dimensionfunctions} (denoted by $F$, and later denoted by $W(A)$ elsewhere).
The difference is that $W(A)$ is constructed using positive elements from matrix algebras over $A$, whereas $\Cu(A)$ uses positive elements from the stabilization of $A$ (and thus $\Cu(A)$ is bigger).
This paper deals primarily with $\Cu(A)$ because it has better structural properties, developed in \cite{CowardElliottIvanescu}; most strikingly, it behaves better with inductive limits.
However, the main results here for $\Cu(A)$ have equally strong analogues for $W(A)$ (in Section \ref{OldCu}), which are derived as corollaries.

\begin{defn}
The Cuntz semigroup of a $C^*$-algebra $A$ is 
\[ \Cu(A) := (A \tens \K)_+ / \cueq. \]
The equivalence class of $a \in (A \tens \K)_+$ is denoted by $[a]$.
This is a semigroup using the addition defined by
\[ [a] + [b] := [a' + b'], \]
where $a \cueq a', b \cueq b'$ and $a' \perp b'$ (such elements can always be found, an the equivalence class of $a' + b'$ does not depend on the particular choice of $a',b'$).
An order is defined on $\Cu(A)$ by
\[ [a] \leq [b] \]
if $a \culeq b$.
\end{defn}

One important property of $\Cu(A)$ is the existence of suprema of increasing sequences:

\begin{prop}(\ccite{Theorem 1 (i)}{CowardElliottIvanescu})
Let $A$ be a $C^*$-algebra.
For any increasing sequence $([c_n]) \in \Cu(A)$, a supremum $[c] \in \Cu(A)$ exists.
\end{prop}

Another useful relation on elements of $\Cu(A)$ is ``compact containment'', given by $[a] \ll [b]$ if whenever $([c_n]) \subset \Cu(A)$ is an increasing sequence and
\[ [b] \leq \sup [c_n], \]
it follows that $[a] \leq [c_n]$ for some $n$.

For a positive element $a$ and a real number $\e > 0$, using the real function $f(t) = (t-\e)_+$, we set
\[ (a-\e)_+ := f(a) \]
using functional calculus.
For $a \in (A \tens \K)_+$, we always have
\[ [(a-\e)_+] \ll [a] \]
(this is essentially contained in \ccite{Theorem 1}{CowardElliottIvanescu}).

In \ccite{Theorem 2}{CowardElliottIvanescu}, it was shown that $\Cu$ is a sequentially continuous functor from the category of $C^*$-algebras to a category of semigroups with certain additional structure.
In particular, if $\phi:A \to B$ is a $*$-homomorphism between $C^*$-algebras then it induces a map between $\Cu(A)$ and $\Cu(B)$, and this map preserves addition, the order $\leq$, and the relation $\ll$.
Additionally, we have the following result about inductive limits.

\begin{prop}(From the proof of \ccite{Theorem 2}{CowardElliottIvanescu}) \alabel{CuIndLimits}
Let
\[ A_1 \labelledrightarrow{\phi_1^2} A_2 \labelledrightarrow{\phi_2^3} \cdots \rightarrow A \]
be an inductive limit of $C^*$-algebras.
Let $[a] \in \Cu(A)$.
\begin{enumerate}
\item[(i)] There exists elements $[a_i] \in \Cu(A_i)$ for each $i$ such that $[\phi_i^{i+1}(a_i)] \leq [a_{i+1}]$ and
\[ [a] = \sup [\phi_i^\infty(a_i)]. \]
\item[(ii)] Let $[a_i]$ be as in (i) and let $[b_i] \in \Cu(A_i)$ for each $i$, such that $[\phi_i^{i+1}(b_i)] \leq [b_{i+1}]$ and 
\[ [b] = \sup [\phi_i^\infty(b_i)]. \]
If $[a] \leq [b]$ then for any $[a'] \ll [a_i]$ in $\Cu(A_i)$, there exists $j \geq i$ such that
\[ [\phi_i^j(a')] \ll [b_j]. \]
\end{enumerate}
\end{prop}

The Murray-von Neumann semigroup is constructed using an order relation on projections; the order relation is defined now.

\begin{defn}
Let $A$ be a $C^*$-algebra.
The order $\mvneq$ is defined on the projections of $A$ by the following.
For projections $p,q \in A$, we have $p \mvneq q$ if there exists $v \in A$ such that $p=v^*v$ and $vv^* = q$.
The pre-order $\mvnleq$ is defined by $p \mvnleq q$ if there exists $q' \leq q$ such that $p \mvneq q'$.
\end{defn}

It turns out that for projections $p,q \in A \tens \K$, we have $p \mvnleq q$ if and only if $p \culeq q$ (by \ccite{Proposition 2.1}{Rordam:UHFII}).
Consequently, if $A$ is stably finite, we have $p \mvneq q$ if and only if $p \cueq q$.

\begin{defn}
The Murray-von Neumann semigroup of a $C^*$-algebra $A$ is defined by
\[ V(A) := \{\text{projections in }A \tens \K\} / \mvneq \]
Addition is given by
\[ [p] + [q] = [p'+q'] \]
where $p \mvneq p', q \mvneq q'$ and $p' \perp q'$.
\end{defn}

When $A$ is stably finite, $V(A)$ is evidently a subsemigroup of $\Cu(A)$.
The following result allows us to easily identify elements of this subsemigroup.

\begin{prop}(\ccite{Theorem 3.5}{BrownCiuperca}) \alabel{BCmagic}
Let $A$ be a stably finite $C^*$-algebra and let $a \in (A \tens \K)_+$.
The following are equivalent:
\begin{enumerate}
\item[(i)] $[a] \ll [a]$ in $\Cu(A)$.
\item[(ii)] $[a] \in V(A)$.
\item[(iii)] The spectrum of $a$ is contained in $\{0\} \cup [\e,\infty)$ for some $\e > 0$.
\end{enumerate}
\end{prop}

It is well-known that, if we have an inductive limit $A = \lim A_i$ of $C^*$-algebras and $[p] \in V(A)$ then $[p]$ is the image of some $[p'] \in V(A_i)$ for some $i$ (this can also be derived, in the stably finite case, from Proposition \ref{BCmagic} (i) and Proposition \ref{CuIndLimits}).

This section is concluded with the following technical result, which allows us to conclude that if $a,b \in C_0(X,A)_+$ satisfy $a|_Y \culeq b|_Y$ for some closed set $Y$, then given $\e > 0$, there exists an open subset $U$ containing $Y$ such that
\[ (a-\e)_+|_{\overline{U}} \culeq b|_{\overline{U}}. \]

\begin{lem}\alabel{CuSemiprojective}
Let $A$ be a $C^*$-algebra and
\[ I_1 \subseteq I_2 \subseteq \cdots \]
be an increasing sequence of ideals of $A$.
Let $a,b \in A_+$ be elements such that
\[ a+\overline{\bigcup_n I_n} \culeq b+\overline{\bigcup_n I_n} \]
in $\Cu(A/\overline{\bigcup_n I_n})$.
Then given $\e > 0$, there exists $n$ such that
\[ (a-\e)_+ + I_n \culeq b+I_n \]
in $\Cu(A/I_n)$.
\end{lem}

\begin{proof}
Let $I = \overline{\bigcup_n I_n}$.
There exists $s+I \in A/I$ such that
\[ \e/2 > \|(s+I)^*(b+ I)(s+I) - (a+I)\| = \|(s^*bs-a)+I\|<\e/2. \]

Hence, for some $n$, we have
\[ \|(s^*bs-a)+I_n\|<\e. \]
By \ccite{Proposition 2.2}{Rordam:UHFII}, it follows that
\[ (a-\e)_+ +I_n \culeq (s^*bs) + I_n \culeq b+I_n, \]
in $A/I_n$.
\end{proof}

\subsection{Approximately subhomogeneous algebras and slow dimension growth}

A $C^*$-algebra is said to be subhomogeneous if there is a finite bound on the dimensions of its irreducible representations.
It is approximately subhomogeneous (ASH) if it can be written as a direct limit of subhomogeneous $C^*$-algebras.
The ASH algebras for which the results of this paper hold are those which are simple, unital, and have slow dimension growth.
The notion of dimension growth for ASH algebras is based on the fact, proven in \cite{NgWinter:ash}, that every unital separable ASH algebra is the direct limit of recursive subhomogeneous (RSH) $C^*$-algebras, as introduced in \ccite{Definitions 1.1 and 1.2}{Phillips:rsh}.
Slow dimension growth for ASH algebras has seen different definitions in the literature (\ccite{Definition 1.1}{Phillips:arsh} and \ccite{Definition 3.2}{Toms:rigidity}); however, for simple, unital algebras, the definitions are known to be equivalent.
We collect these known equivalences in the following proposition; conditions (iv) and (v) are what will be used in this paper.

\begin{prop}\alabel{ASH-SlowCharacterization}
Let $A$ be a simple, unital, non-type I ASH algebra.
The following are equivalent.
\begin{enumerate}
\item[(i)]
$A$ has slow dimension growth in the sense of \ccite{Definition 1.1}{Phillips:arsh};
\item[(ii)]
$A$ has slow dimension growth in the sense of \ccite{Definition 3.2}{Toms:rigidity};
\item[(iii)]
$A$ can be written as a direct limit of a system $(A_i,\phi_i^{i+1})$ that has slow dimension growth as in \ccite{Definition 3.2}{Toms:rigidity}, and additionally, the connecting maps $\phi_i^{i+1}$ are unital and injective.
\item[(iv)]
$A$ can be written as a direct limit of RSH algebras $(A_i,\phi_i^{i+1})$, such that for every $i$, every non-zero $c \in A_i$, and every $N$, there exists $j \geq i$ such that, for every $\omega$ in the total space of $A_j$,
\[ \rank \ev_\omega(\phi_i^j(c)) \geq Nd(\omega) \]
(where the evaluation map $\ev_\omega$ and the topological dimension function $d(\cdot)$ are as defined in \ccite{Definition 1.2}{Phillips:rsh}).
\item[(v)]
$A \iso A \tens \jsZ$.
\end{enumerate}
\end{prop}

\begin{proof}
(i) $\Rightarrow$ (ii) is easy to see, by setting $p$ to be the unit in \ccite{Definition 1.1}{Phillips:arsh}.
(ii) $\Rightarrow$ (v) follows from \ccite{Theorem 1.2}{Toms:rigidity} and \ccite{Corollary 6.4}{Winter:Perfect}.
(v) $\Rightarrow$ (iii) is shown by (the proof of) \ccite{Theorem 5.5}{TomsWinter:ssa}.
(iv) $\Rightarrow$ (i) is obvious.

Let us prove (iii) $\Rightarrow$ (iv).
Let $(A_i,\phi_i)$ be the system given by (iii), and let us be given a non-zero element $c \in A_i$ and some $N$.
Since $c$ and $c^*c$ have the same rank at each point, we may assume that $c \geq 0$.

By \ccite{Lemma 1.5}{Phillips:arsh}, we may find $j_1 \geq i$ such that, for every $\omega$ in the total space of $A_{j_1}$, we have
\[ \ev_\omega(\phi_i^{j_1}(c)) \neq 0. \]
By \ccite{Theorem 4.6}{Toms:comparison}, we can find $M > 0$ such that $1_{A_{j_1}} \culeq c^{\dsum M}$.
Since the system $(A_i,\phi_i)$ satisfies \ccite{Definition 3.2}{Toms:rigidity}, we can find $j \geq j_1$ such that
\[ \rank(\ev_\omega(1_{A_j})) \geq NMd(\omega) \]
for all $\omega$ in the total space of $A_j$.
Since $1_{A_j} \culeq \phi_i^j(c)^{\dsum M}$, we must have
\[ NMd(\omega) \leq M\rank \ev_\omega(\phi_i^j(c)). \]
\end{proof}
\section{General considerations\alabel{GeneralConsid}}

\begin{prop} \alabel{Phi-ll-lsc}
Let $A$ be a $C^*$-algebra, let $X$ be a locally compact Hausdorff space, and let $[a] \in \Cu(C_0(X,A))$.
Then the map $x \mapsto [a(x)]$ is lower semicontinuous with respect to $\ll$; that is, for every $[b] \in \Cu(A)$, the set
\[ \{x \in X: [b] \ll [a(x)]\} \]
is open.
\end{prop}

\begin{proof}
Let $x \in X$ be such that $[b] \ll [a(x)]$.
Then, for some $\e > 0$ we have $[b] \ll [(a-\e)_+(x)]$.
Let $U$ be a neighbourhood of $x$ such that for $y \in U$, $\|a(y)-a(x)\| < \epsilon$.
Then by \ccite{Proposition 2.2}{Rordam:UHFII}, we have 
\[ (a(x)-\epsilon)_+ \culeq a(y), \]
for all $y \in U$.
Thus, $U$ is an open neighbourhood of $x$ contained in $\{x \in X: [b] \ll [a(x)]\}$.
\end{proof}

\begin{rmk}\alabel{LLSemicontinuity}
It follows that, also, the set
\[ \{x \in X: [b] \ll [a(x)], [b] \neq [a(x)]\} \]
is open.
To see this, note that this set is the union of the sets
\[ \{x \in X: [c] \ll [a(x)]\} \]
given by all $[c]$ satisfying $[b] < [c]$.
\end{rmk}

\begin{prop}\alabel{PointwiseProj}
Let $A$ be a separable stably finite $C^*$-algebra and $X$ a compact Hausdorff space.
Suppose that $[a] \in \Cu(C(X,A))$ is such that $[a(x)]$ is equal to the same $[p] \in V(A)$ for all $x \in X$
Then $[a] \in V(C(X,A))$.
\end{prop}

\begin{proof}
We have that $\chii_{(0,\infty)}(a(x)) \in A \tens \K$ is defined at all $x \in X$; we must show that it is continuous.
For $x \in X$, $a(x)$ has spectrum contained in $\{0\} \cup [\delta,\infty)$ for some $\delta > 0$ (by Proposition \ref{BCmagic})
So, there exists $f \in C_0((0,\infty))$ such that $p := \chii_{(0,\infty)}(a(x)) = f(a(x))$.
Given $\eta > 0$, there exists a neighbourhood $U$ of $x$ such that $\|f(a(x)) - f(a(y))\| < \eta$ for $y \in U$.
Thus,
\[ \|p - f(a(y))pf(a(y))\| < 3\eta(1+\eta). \]
For a given $\epsilon > 0$, if $\eta$ is sufficiently small (how small depends only on $\epsilon$), this implies that there is a projection $q \in \hered{f(a(y))} \subseteq \hered{a}$ with $\|p - q\| < \epsilon$.
In this case (assuming $\epsilon \leq 1$), we have $p \mvneq q \leq \chii_{(0,\infty)}(a(y))$.
But, since $p \mvneq \chii_{(0,\infty)}(a(y))$ and $A \tens \K$ is finite, we must have that $q = \chii_{(0,\infty)}(a(y))$.
This shows that, for $y \in U$, $\chii_{(0,\infty)}(a(y))$ is distance at most $\epsilon$ from $p = \chii_{(0,\infty)}(a(x))$, thus showing continuity of $x \mapsto \chii_{(0,\infty)}(a(x))$.
\end{proof}

In the next result, we no longer assume that $X$ is compact.
Part (ii) of this result was essentially obtained in \ccite{Lemma 5.2}{CommutativeCuntz}, for the case that $A = \K$.

\begin{prop} \alabel{CuntzForC0Projs}
Let $A$ be a separable stably finite $C^*$-algebra and $X$ a second countable locally compact Hausdorff space.
Suppose that $[a] \in \Cu(C_0(X,A))$ is such that $[a(x)] \in V(A)$ for all $x \in X$ and the map $x \mapsto [a(x)]$ is continuous (ie.\ locally constant).
Set $p(x) = \chii_{(0,\infty)}(a(x))$ for each $x$.
Then:
\begin{enumerate}
\item[(i)] $p$ is a continuous map from $X$ to the projections in $A \tens \K$,
\item[(ii)] If $f \in C_0(X)_+$ is strictly positive then $fp \in C_0(X,A \tens \K)_+$ and
\[ [a] = [fp] \]
in $\Cu(C_0(X,A))$.
\item[(iii)] If $q$ is another continuous map from $X$ to the projections in $A \tens \K$ and if $g \in C_0(X)_+$ is strictly positive then $gq \in C_0(X, A \tens \K)_+$ also, and we have
\[ [fp] \leq [gq] \]
if and only if for every compact set $K \subseteq X$, $p|_K \mvnleq q|_K$.
\item[(iv)] Given $q$ and $g$ as in (iii), we have that 
\[ [fp] = [gq] \]
if and only if $[fp] \leq [gq]$ and $p(x) \mvneq q(x)$ for each $x \in X$.
\end{enumerate}
\end{prop}

\begin{proof}
(i): Fix a point $x \in X$; we may find a compact neighbourhood $K$ of $x$.
By the proof of Proposition \ref{PointwiseProj}, we see that $p|_K$ is continuous with values in $A \tens \K$.

(ii) and (iii): It is clear that $fp,gq \in C_0(X,A \tens \K)$.
To show the rest of (ii) and (iii), we shall show more generally that if $[b]$ is another element of $\Cu(C_0(X,A))$ satisfying the condition that the map $x \mapsto [b(x)]$ is locally constant and $q(x) = \chii_{(0,\infty)}(x)$ then $[a] \leq [b]$ if and only if for every compact set $K \subseteq X$, $p|_K \mvnleq q|_K$.

On the one hand, if $[a] \leq [b]$ then $[a|_K] \leq [b|_K]$.
But we have $[a|_K] = [p|_K]$ and $[b|_K] = [q|_K]$, so $p|_K \culeq q|_K$.
Since $p|_K, q|_K$ are projections in $C(K,A \tens \K)$, it must be the case that $p|_K \mvnleq q|_K$.

Conversely, suppose that $p|_K \mvnleq q|_K$ for each compact subset $K$ of $X$.
Given $\e > 0$, let $K$ be the support of $(a-\e/2)_+$, which is compact.
Let $v \in C(K,A \tens \K)$ be such that $p|_K = v^*v, vv^* = q|_K$.
Using Proposition \ref{BCmagic} as in the proof of Proposition \ref{PointwiseProj}, let $f$ be a continuous function such that $f(b)(x) = q(x)$ for $x \in K$.
Let $h \in C_0(X)_+$ be such that $h$ is zero outside of $K$ and $h$ takes the value $1$ on the support of $(a-\e)_+$; this is possible since the support of $(a-\e)_+$ is compactly contained in that of $(a-\e/2)_+$.
Then we have
\begin{align*}
(a-\e)_+ = (a-\e)_+^{1/2}hph(a-\e)_+^{1/2} &= (a-\e)_+^{1/2}hv^*f(b)vh(a-\e)_+^{1/2} \\
\culeq f(b) \culeq b.
\end{align*}
Since $\e$ is arbitrary, this shows that $a \culeq b$.

(iv) It is clear that if $[fp] = [gq]$ then $[fp] \leq [gq]$ and $p(x) \mvneq q(x)$ for all $x$.
Conversely, for each $K \subseteq X$ compact, we have a partial isometry $v \in C(K,\K)$ such that $p|_K = v^*v, vv^* \leq q|_K$.
If $vv^* \neq q|_K$ then for some $x \in K$ we have $vv^*(x) < q(x)$, yet $q(x) \mvneq p(x) \mvneq vv^*(x)$, contradicting stable finiteness of $A$.
Hence, $vv^* = q|_K$.
It follows from (iii) that $[p] = [q]$.
\end{proof}

Recall from the introduction that
\[ V_c(X,A) := \{\text{projection valued functions } p \in C_b(X, A \tens \K)\}/ \sim_c \]
where $p\sim_c q$ if $p|_K \mvneq q|_K$ for every compact $K$.
We denote by $\langle p \rangle$ the equivalence class of the projection-valued function $p \in C_b(X,A \tens \K)$, and more generally, if $a:X \to A \tens \K$ is such that $[a(x)] \in V(A)$ for each $x$ and the induced map $x \to [a(x)]$ is locally constant, then
\[ \langle a \rangle := \langle \chii_{(0,\infty)} a \rangle, \]
where the functional calculus is done pointwise.

\begin{cor}\alabel{VcDescr}
Let $A$ be a separable stably finite $C^*$-algebra and $X$ a $\sigma$-compact, locally compact Hausdorff space.
The semigroup $V_c(X,A)$ may be identified with each of the following:
\begin{enumerate}
\item[(i)] The subsemigroup of $\Cu(C_0(X,A))$ consisting of all $[a]$ for which the map $x \mapsto [a(x)]$ is continuous with range in $V(A)$.
\item[(ii)] The inverse limit
\[ \varprojlim_{K \mathrm{ compact}, K \nearrow X} V(C(K,A)) \]
where the connecting maps are given by restriction, $V(C(K,A)) \to V(C(L,A)): p \mapsto p|_L$ when $L \subseteq K$.
\end{enumerate}
\end{cor}

\begin{proof}
Proposition \ref{CuntzForC0Projs} shows that $[a] \mapsto \langle \chii_{(0,\infty)}(a) \rangle$ (functional calculus taken pointwise) is a well-defined injective map from the subsemigroup of $\Cu(C_0(X,A))$ described in (i).
Moreover, since $C_0(X,A)_+$ has a strictly positive element, the map is clearly onto.

It is clear, from the definition of $\sim_c$, that
\[ \langle p \rangle \mapsto ([p|_K] \in V(C(K,A)))_{K \subset X \text{ compact}} \]
is a well-defined, injective map from $V_c(X,A)$ to $\varprojlim V(C(K,A))$.
To see that it is surjective, let $([q_K]) \in \varprojlim V(C(K,A))$; that is to say, we have some element $[q_K]$ of $V(K,A)$ for each compact subset $K$ of $X$, and if $K \subseteq L$ then $q_L|_K \mvneq q_K$.
Write $X$ as a union of a sequence $(K_i)_{i=1}^\infty$ of compact sets, such that for each $i$, $K_i$ is contained in the interior of $K_{i+1}$.
Then, we may find $f_i \in C_0(X)_+$ such that $f_i|_{K_{i+1}^c} = 0$ and $f_i$ is strictly positive on $K_i$, and set
\[ [a_i] = [f_iq_{K_{i+1}}] \in \Cu(C_0(X,A)). \]
One can easily verify that $[a_i] \leq [a_{i+1}]$ for each $i$, and if
\[ [a] = \sup [a_i] \]
then $[a|_{K_i}] = [q_{K_i}]$ for each $i$.
Since each compact subset $K$ of $X$ is contained in some $K_i$, it follows that $[a|_K] = [q_K]$ for each compact $K$.
\end{proof}

If $Y$ is a subset of $X$, arising as the intersection of a closed and an open subset of $X$, then there is an induced map $\Cu(C_0(X,A)) \to \Cu(C_0(Y,A))$ (as described in \ccite{Section 2.1, page 4}{CommutativeCuntz}, see also \cite{CiupercaRobertSantiago}).
The map takes $V_c(X,A)$ into $V_c(Y,A)$, by
\[ \langle a \rangle \to \langle a|_Y \rangle. \]

\section{Cuntz comparison in $C_0(X,ASH)$\alabel{ComparisonSection}}

The following result determines Cuntz comparison for a separable $C^*$-algebra arising as a commutative algebra tensored with a simple ASH algebra with slow dimension growth.

\begin{thm}\alabel{ASH-Comparison}
Let $A$ be a simple, unital, $\jsZ$-stable ASH algebra and let $X$ be a second countable, locally compact Hausdorff space.
Let $a,b \in C_0(X,A\tens \K)_+$ be positive elements.
Then $a \culeq b$ if and only if
\begin{enumerate}
\item[(i)] for each $x \in X$, $a(x) \culeq b(x)$, and
\item[(ii)] for each $[p] \in V(A)$, we have
\[ \langle a|_{\{x: [a(x)]=[p]=[b(x)]\}} \rangle = \langle b|_{\{x: [a(x)]=[p]=[b(x)]\}} \rangle \]
in $V_c(\{x: [a(x)]=[p]=[b(x)]\},A)$.
\end{enumerate}
\end{thm}

It should be noted that the set $\{x \in X: [a(x)]=[p]=[b(x)]\}$ appearing in condition (ii) is locally compact, as it can be written as a difference of open sets:
\[ \{x \in X: [a(x)]=[p]=[b(x)]\} = \{x \in X: [p] \ll [a(x)]\} \setminus \{x \in X: [p] \ll [b(x)], [p] \neq [b(x)]\}; \]
the latter set is open by Remark \ref{LLSemicontinuity}.

Theorem \ref{ASH-Comparison} will by proven by reducing to the following result, which provides a partial determination of Cuntz comparison for a commutative algebra tensored with an RSH algebra.
We once again refer to \ccite{Definition 1.1 and 1.2}{Phillips:rsh} for the notation associated to an RSH algebra.

\begin{lem}\alabel{RSH-Embedding}
Let $X$ be a finite dimensional locally compact Hausdorff space and let $Y \subseteq X$ be a closed subset.
Let $R$ be an RSH algebra with finite dimensional total space $\Omega$, and let $\sigma:R \to C(\Omega,\K)$ be the canonical representation of $R$.
Suppose that $a,b \in C_0(X,R \tens \K)_+$ satisfy
\begin{enumerate}
\item[(i)] for all $x \in X \setminus Y$ and all $\omega \in \Omega$,
\[ \rank \sigma(a(x))(\omega) + \frac{\dim X + d(\omega) - 1}{2} \leq \rank \sigma(b(x))(\omega); \text{ and} \]
\item[(ii)] there exists $s \in C_0(Y,R \tens \K)$ such that $s^*s = a|_Y$ and $ss^* \in \hered{b|_Y}$.
\end{enumerate}
Then there exists $\tilde{s} \in C_0(X, R \tens \K)$ such that $\tilde{s}|_Y = s$, $\tilde{s}^*\tilde{s} = a$ and $\tilde{s}\tilde{s}^* \in \hered{b}$.
\end{lem}

\begin{proof}
Fix a decomposition of the RSH algebra $R$ (see \ccite{Definition 1.1}{Phillips:rsh}).
This result shall be proven inductively on the length of decomposition of $R$.
If $R$ is simply a matrix algebra, then the result is exactly \ccite{Corollary 3.3}{CommutativeCuntz}.
When $R$ has length $0$, this means that $R = C(\Omega,M_n)$ for some $n$.
Using the case that $R$ is a matrix algebra as the base case, the same inductive step handles the case of length $0$ and greater lengths.

The inductive step is as follows.
Let $R$ be given by the pull-back
\begin{equation}
\begin{array}{rcl}
R & \rightarrow & C(\Gamma, M_n) \\
{\scriptstyle \lambda} \downarrow && \quad\ \downarrow {\scriptstyle f \mapsto f|_{\Gamma_0}} \\
R' & \labelledrightarrow{\rho} & C(\Gamma_0, M_n),
\end{array}
\alabel{RSH-Embedding-Diagram}
\end{equation}
where $R'$ is an RSH algebra whose length is one less than that of $R$, and $\rho$ is a unital $*$-homomorphism (or, when the length of $R$ is $0$, then $R' = M_n$, $\Gamma_0$ is any singleton subset of $\Gamma$, and $\rho$ is the obvious identification of $M_n$ with $C(\Gamma_0,M_n)$).
Abusing notation slightly, we will denote by $\sigma$ the map from $R$ to $C(\Gamma, M_n)$ in the top row of \eqref{RSH-Embedding-Diagram}.

By induction, there exists $\tilde{s}_0 \in C_0(X,R')$ such that $\tilde{s}_0|_Y = (\id_{C_0(Y)} \tens \lambda)(s|_Y)$, $\tilde{s}_0^*\tilde{s}_0 = (\id_{C_0(X)} \tens \lambda)(a)$ and $\tilde{s}_0\tilde{s}_0^* \in \hered{(\id_{C_0(X)} \tens \lambda)(b)}$.
Let $t:(X \times \Gamma_0) \cup (Y \times \Gamma) \to M_n \tens \K$ be defined by
\begin{align*}
t|_{X \times \Gamma_0} &= \rho(\tilde{s}_0), \\
t|_{Y \times \Gamma} &= \sigma(s).
\end{align*}
By the commutativity of \eqref{RSH-Embedding-Diagram}, we see that $t$ is well-defined on $Y \times \Gamma_0$ (and therefore, continuous on its entire domain).
Applying \ccite{Corollary 3.3}{CommutativeCuntz}, we may extend $t$ to $\tilde{t} \in C_0(X \times \Gamma, M_n \tens \K)$ such that $\tilde{t}^*\tilde{t} = \sigma(a)$ and $\tilde{t}\tilde{t}^* \in \hered{\sigma(b)}$.
It follows that $\tilde{s} := (\tilde{s}_0, \tilde{t})$ is an element of $C_0(X, R \tens \K)$ and that it satisfies $\tilde{s}|_Y = s$, $\tilde{s}^*\tilde{s} = a$ and $\tilde{s}\tilde{s}^* \in \hered{b}$.
\end{proof}

In order to apply Lemma \ref{RSH-Embedding} in the proof of Theorem \ref{ASH-Comparison}, we shall use the following result, which allows us to move from a simple, non-type I limit of RSH algebras to one of the building blocks.

\begin{lem}\alabel{ConditionToBuildingBlocks}
Let
\[ X_1 \labelledleftarrow{\alpha_2^1} X_2 \labelledleftarrow{\alpha_3^2} \cdots \leftarrow X \]
be a projective limit of locally compact Hausdorff spaces, and let $A$ be a stably finite $C^*$-algebra given as an inductive limit
\[ A_1 \labelledrightarrow{\phi_1^2} A_2 \labelledrightarrow{\phi_2^3} \cdots \rightarrow A. \]
Assume that each map $\alpha_{i+1}^i$ is proper and surjective while each map $\phi_i^j$ is injective.
Use $\Phi_i^j:C_0(X_i,A_i) \to C_0(X_j,A_j)$ to denote the map $f \mapsto \phi_i^j \circ f \circ \alpha_j^i$, so that we have an inductive limit decomposition 
\[ C_0(X_1,A_1) \labelledrightarrow{\Phi_1^2} C_0(X_2,A_2) \labelledrightarrow{\Phi_2^3} \cdots \rightarrow C_0(X,A). \]
Let $a \in C_0(X_1,A_1)_+$.
Let $b_i \in C_0(X_i,A_i)_+$ be such that $\Phi_i^{i+1}(b_i) \culeq b_{i+1}$ for each $i$ and set
\[ [b] = \sup_i [\Phi_i^\infty(b_i)] \in \Cu(C_0(X,A)). \]
Suppose that, for all $x \in X$, we have $[\Phi_1^\infty(a)(x)] \leq [b(x)]$.
Let $\e > 0$.
\begin{enumerate}
\item[(i)] If we have, for all $[p] \in V(A)$ that
\[ \langle \Phi_1^\infty(a)|_{\{x: [\Phi_1^\infty(a)(x)] = [p] = [b(x)]\}} \rangle = \langle b|_{\{x: [\Phi_1^\infty(a)(x)] = [p] = [b(x)]\}} \rangle \]
in $V_c(\{x: [\Phi_1^\infty(a)(x)] = [p] = [b(x)]\}, A)$ then, given $\e > 0$, there exists $i \geq 1$ and a closed set $Y$ in $X_i$ such that,
\[ [\Phi_1^i((a-\e)_+)|_Y] \leq [b_i|_Y] \]
in $\Cu(C_0(Y,A_i))$ and, for $x \in (\alpha_\infty^i)^{-1}(X_i \setminus Y^\circ)$,
\[ [\Phi_1^\infty(a)(x)] < [b(x)] \]
in $\Cu(A)$.

\item[(ii)]
If we in fact have $[\Phi_1^\infty(a)(x)] < [b(x)]$ for all $x \in X$ then there exists $i$ such that, for all $x \in X_i$, either $\Phi_1^i((a-\e)_+)(x)=0$ or
\[ [\Phi_1^i((a-\e)_+)(x)] < [b_i(x)]. \]

\item[(iii)]
Suppose that the system $(A_i,\phi_i^{i+1})$ is as in Proposition \ref{ASH-SlowCharacterization} (iv).
Suppose that for all $x \in X$, we in fact have $[\Phi_1^\infty(a)(x)] < [b(x)]$.
Then given $N > 0$, there exists $i$ such that, for all $x \in X$, either $\Phi_1^i((a-\e)_+)(x)=0$ or for all $\omega$ in the total space of $A_i$,
\[ \rank\, \ev_\omega(\Phi_1^i((a-\e)_+(x))) + Nd(\omega) \leq \rank\, \ev_\omega(b_i(x)). \]
\end{enumerate}
\end{lem}

\begin{proof}
(i)
Let $\eta = \e/2$.
Let $Y_0 := \{x \in X: [\Phi_1^\infty((a-\eta)_+)(x)] = [b(x)] \neq 0\}$.

\begin{claim}
The set $Y_0$ is compact.
\end{claim}

\begin{proof}[Proof of claim.]
For $y \in \overline{Y_0}$, we can for any $\dl > 0$ find $x \in Y_0$ such that $\|b(x)-b(y)\|<\dl$ and $\|\Phi_1^\infty(a)(x)-\Phi_1^\infty(a)(y)\|<\eta$.
So, using \ccite{Proposition 2.2}{Rordam:UHFII}, we have
\begin{equation}
[(b(y)-\dl)_+] \leq [b(x)] = [\Phi_1^\infty((a-\eta)_+)(x)] \leq [\Phi_1^\infty(a)(y)] \leq [b(y)].
\alabel{ConditionToBuildingBlocks-abeqn}
\end{equation}
Taking $\dl$ arbitrarily small, this shows that $[b(y)]=[\Phi_1^\infty(a)(y)] \neq 0$, ie.\ $y \in Y_0$, and so $Y_0$ is closed.

Since we know that $(a-\eta)_+$ has compact support, and this support contains $Y_0$ by definition, it follows that $Y_0$ is compact.
\end{proof}

Evident from \eqref{ConditionToBuildingBlocks-abeqn} is the fact that, if $x,y \in Y_0$ are such that $\|\Phi_1^\infty(a)(x)-\Phi_1^\infty(a)(y)\|<\eta$ then $[b(x)]=[b(y)]$.
Consequently, $Y_0$ is the disjoint union of relatively open sets of the form
\[ Y_0 \cap \{x: [b(x)]=[b(y)]=[p]\}, \]
(where $[p]$ varies over $V(A)$).
Since $Y_0$ is compact, it is in fact a finite union of such sets, and as such, we see that
\[ [\Phi_1^\infty(a)|_{Y_0}] = [b|_{Y_0}] \]
from the hypothesis.

Since ${Y_0}$ is compact, we have $[b|_{Y_0}] \ll [b|_{Y_0}]$.
By Proposition \ref{CuIndLimits} (ii), it follows that, for some $i$ we have
\[ [\Phi_1^i(a)|_{\alpha_\infty^i({Y_0})}] = [b_i|_{\alpha_\infty^i({Y_0})}]. \]

Now, let $Y_1 := \{x \in X: \phi_i^\infty(a_i(x)) = 0\}$, so that $Y_1$ is closed, and in fact, $Y_1 \cup \{\infty\}$ is closed in the one-point compactification of $X$.
From this it follows that $\alpha_\infty^i(Y_1)$ is closed, and it is clear that
\[ [\Phi_1^i(a)|_{\alpha_\infty^i(Y_0 \cup Y_1)}] = [b_i|_{\alpha_\infty^i(Y_0 \cup Y_1)}]. \]
By Lemma \ref{CuSemiprojective}, there exists a closed set $Y$ in $X$, such that $Y_0 \cup Y_1$ is contained in the interior of $Y$, and
\[ [\Phi_i^{j_1}((a_i-\e)_+)|_Y] \leq [b_i|_Y]. \]

By our choice of $Y_0$ and $Y_1$, we clearly have $[\Phi_1^\infty(a)(x)] < [b(x)]$ for all $x \in (\alpha_\infty^i)^{-1}(X_i \setminus Y^\circ)$.

(ii)
Let $\eta = \e/2$.
For $x \in X$, we have $[\Phi_1^\infty((a-\eta)_+(x))] \ll [b(x)]$ and $\Phi_1^\infty((a-\eta)_+(x))] \neq [b(x)]$, and therefore by Proposition \ref{CuIndLimits} (ii) (and considering separately the cases that $[\Phi_1^\infty(a)(x)]$ is compact or not), we have for some $i_x$ that
\begin{align*}
[\phi_1^{i_x}((a-\eta)_+(\alpha_\infty^{1}(x)))] &\ll [b_{i_x}(\alpha_\infty^{i_x}(x))], \text{ and} \\
[\phi_1^{i_x}((a-\eta)_+(\alpha_\infty^{1}(x)))] &\neq [b_{i_x}(\alpha_\infty^{i_x}(x))].
\end{align*} 
Hence, for some $\dl>0$,
\[ [\phi_1^{i_x}((a-\eta)_+(\alpha_\infty^{1}(x)))] < [(b_{i_x}-\dl)_+\alpha_\infty^{i_x}(x))]. \]

By \ccite{Proposition 2.2}{Rordam:UHFII}, let $U_x$ be an open neighbourhood of $\alpha_\infty^{i_x}(x)$ in $X_{i_x}$ such that, for all $y \in U_x$ we have
\[ \phi_1^{i_x}((a-\e)_+(\alpha_{i_x}^1(y))) \culeq \phi_1^{i_x}((a-\eta)_+(\alpha_\infty^{1}(x))) \]
and
\[ (b_{i_x}-\dl)_+(\alpha_\infty^{{i_x}}(x)) \culeq b_{i_x}(y). \]
Thus, for $y \in U_x$, we have
\[ [\phi_1^{i_x}((a-\e)_+(y))] < [b_{i_x}(y)]. \]

The sets $(\alpha_\infty^{i_x})^{-1}(U_x)$ form an open cover of the support of $\Phi_1^\infty((a-\e)_+)$, which is a compact set, whence there are $x_1,\dots,x_\ell \in X$ such that the sets $(\alpha_\infty^{i_{x_1}})^{-1}(U_{x_1}), \dots, (\alpha_\infty^{i_{x_\ell}})^{-1}(U_{x_\ell})$ cover the support of $\Phi_1^\infty((a-\e)_+)$.
Letting $i$ be the maximum of $i_{x_1},\dots,i_{x_\ell}$, (ii) follows since, for any $x \in X_i$, if $\Phi_1^{i}((a-\e)_+)(x) \neq 0$ then $x$ is contained in $(\alpha_i^{i_{x_1}})^{-1}(U_{x_1}) \cup \cdots \cup (\alpha_i^{i_{x_\ell}})^{-1}(U_{x_\ell})$.

The proof of (iii) runs along similar lines to the proof of (ii).
Letting $\eta = \e/3$, we can first find $j \geq i$ and $\dl > 0$ such that
\[ [\phi_1^j((a-\eta)_+(\alpha_\infty^{1}(x)))] < [(b_{j}-\dl)_+\alpha_\infty^{j}(x))]. \]

\begin{claim}
There exists a non-zero element $c \in (A_j \tens \K)_+$ such that
\[ \phi_1^{j}((a-2\eta)_+)(\alpha_\infty^{1}(x))) \dsum c \culeq (b_j-\delta)_+(\alpha_\infty^{j}(x)). \]
\end{claim}

\begin{proof}[Proof of claim.]
Set $a' :=\phi_1^{j}(a-\eta)_+(\alpha_\infty^{1}(x)))$ and $b'=(b_j-\delta)_+(\alpha_\infty^{j}(x))$, so we have $[a'] < [b']$ and we want to show that $[(a'-\eta)_+] + [c] \leq [b']$ for some non-zero $c$.

Certainly, if $[a'] \in V(A)$ then this follows by \ccite{Proposition 2.2}{PereraToms:recasting}.
Otherwise, $f \in C_0((0,\eta))_+$ be strictly positive, so by Proposition \ref{BCmagic}, 
\[ f(a') \neq 0. \]
However, $f(a') \perp (a'-\eta)_+$ and so
\[ (a-\eta)_+ \dsum f(a') \culeq a' \culeq b'. \]
\end{proof}

By Proposition \ref{ASH-SlowCharacterization} (iv), we can find $i_x \geq j$ such that $\rank\, \ev_\omega(\phi_j^{i_x}(c)) \geq Nd(\omega)$ for all $\omega$ in the total space of $A_{i_x}$, so that we have
\[ \rank\, \ev_\omega \phi_1^{i_x}((a-2\eta)_+)(\alpha_\infty^{1}(x))) + Nd(\omega)  \leq \rank\, \ev_\omega \phi_j^{i_x}(b_j-\delta)_+(\alpha_\infty^{j}(x)). \]
Then, using \ccite{Proposition 2.2}{Rordam:UHFII} as in the proof of (ii), we can obtain a neighbourhood $U_x$ of $\alpha_\infty^{i_x}(x)$ in $X_{i_x}$ such that, for $y \in U_x$,
\[ \rank\, \ev_\omega \Phi_1^{i_x}((a-\e)_+)(y) + d(\omega) \leq \rank\, \ev_\omega b_{i_x}(y). \]

The rest of the proof of (iii) follows a compactness argument as in the proof of (ii).
\end{proof}

\begin{proof}[Proof of Theorem \ref{ASH-Comparison}.]
Let us first prove the theorem in the case that $X$ is finite-dimensional.
Let
\[ A_1 \labelledrightarrow{\phi_1^2} A_2 \labelledrightarrow{\phi_2^3} \cdots \rightarrow A \]
be an inductive system as in Proposition \ref{ASH-SlowCharacterization} (iv).
Set $\Phi_i^j:= \id_{C_0(X)} \tens \phi_i^j:C_0(X,A_i) \to C_0(X,A_j)$.

By Proposition \ref{CuIndLimits} (i), we may find $[a_i], [b_i] \in \Cu(C_0(X,A_i))$ such that
\[ [\Phi_i^{i+1}(a_i)] \leq [a_{i+1}],\ [\Phi_i^{i+1}(b_i)] \leq [b_{i+1}] \]
in $\Cu(C_0(X,A_{i+1}))$ and
\[ [a] = \sup [\Phi_i^\infty(a_i)],\ [b] = \sup [\Phi_i^\infty(b_i)] \]
in $\Cu(C_0(X,A))$.
Given $i \in \mathbb{N}$ and $\e > 0$, let us show that
\[ [\Phi_i^\infty((a_i-\e)_+)] \leq [b]. \]

Let $\eta = \e/2$.
By Lemma \ref{ConditionToBuildingBlocks} (i), we may find $j_1 \geq i$ and a closed set $Y$ in $X$ such that
\[ [\Phi_i^\infty(a_i)(x)] < [b(x)]. \]
for $x \in X \setminus Y^\circ$ and
\[ [\Phi_i^{j_1}((a_i-\eta)_+)|_Y] \leq [b_i|_Y]. \]
Therefore, there exists $s \in C_0(Y,A_i \tens \K)_+$ such that
\[ \Phi_i^{j_1}((a_i-\e)_+) = s^*s \text{ and } ss^* \in \hered{b_i|_Y}. \]

By applying \ref{ConditionToBuildingBlocks} (iii) to $a_i|_{X \setminus Y^\circ}$, we can find $j_2 \geq j_1$ such that, for every $x \in X_{j_2}$ and every $\omega$ in the total space of $A_{j_2}$,
\[ \rank\, \ev_\omega(\phi_i^{j_2}((a_i-\e)_+(x)) + \frac{\dim X + d(\omega) - 1}2 \leq \rank\, \ev_\omega b_i(x).  \]

We may now apply Lemma \ref{RSH-Embedding}, to obtain $\tilde{s} \in C_0(X,A_{j_2} \tens \K)$ such that
\[ \tilde{s}^*\tilde{s} = \Phi_i^{j_2}((a_i-\e)_+) \text{ and } \tilde{s}\tilde{s}^* \in \hered{b_i}. \]
It follows that $[\Phi_i^\infty((a_i-\e)_+)] \leq [b]$, as required.
Since $i$ and $\e$ are arbitrary, this completes the proof in the case that $X$ is finite-dimensional.

For general $X$, we may write $X$ as an inverse limit of finite-dimensional locally compact Hausdorff spaces
\[ X_1 \labelledleftarrow{\alpha_2^1} X_2 \labelledleftarrow{\alpha_3^2} \cdots \leftarrow X, \]
where the connecting maps, $\alpha_j^i$, are proper and surjective.
Use $\Phi_i^j:C_0(X_i,A) \to C_0(X_j,A)$ to denote the map $f \mapsto f \circ \alpha_j^i$.

Again, use Proposition \ref{CuIndLimits} (i) to find $[a_i], [b_i] \in \Cu(C_0(X_i,A))$ such that
\[ [\Phi_i^{i+1}(a_i)] \leq [a_{i+1}],\ [\Phi_i^{i+1}(b_i)] \leq [b_{i+1}] \]
in $\Cu(C_0(X_i,A))$ and
\[ [a] = \sup [\Phi_i^\infty(a_i)],\ [b] = \sup [\Phi_i^\infty(b_i)] \]
in $\Cu(C_0(X,A))$.

Let $i$ and $\e > 0$ be given.
By Lemma \ref{ConditionToBuildingBlocks} (i), we may find $j_1 \geq i$ and a closed set $Y$ in $X_{j_1}$ such that
\[ [\Phi_i^{j_1}((a_i-\e)_+)|_Y] \leq [b_{j_1}|_Y] \]
in $\Cu(C_0(Y,A))$ and, for $x \in (\alpha_\infty^i)^{-1}(X_i \setminus Y^\circ)$,
\[ [\Phi_1^\infty(a_i)(x)] < [b(x)]. \]
in $\Cu(A)$.

Applying Lemma \ref{ConditionToBuildingBlocks} (ii) to $X_i \setminus Y^\circ$, we can find $j_2 \geq j_1$ such that, for $x \in (\alpha_{j_2}^{j_1})^{-1}(X_i \setminus Y^\circ)$,
\[ [\Phi_i^{j_2}((a_i-\e)_+)(x)] < [b_{j_2}(x)]. \]

Then by applying the finite-dimensional case just proven to $\Phi_i^{j_2}((a_i-\e)_+)$ and $b_{j_2}$ in $C_0(X_{j_2},A)$, we see that
\[ [\Phi_i^{j_2}((a_i-\e)_+)] \leq [b_{j_2}], \]
and so
\[ [\Phi_i^\infty((a_i-\e)_+)] \leq [b]. \]
\end{proof}
\section{The Cuntz classes that arise\alabel{DataAttained}}
In Theorem \ref{ASH-Comparison}, it was shown, for a simple ASH algebra $A$ with slow dimension growth, that the Cuntz class of a positive element $a$ in $C_0(X,A \tens \K)_+$ is determined by the value of a certain Cuntz equivalence invariant, consisting of the lower semicontinuous function
\[ x \mapsto [a(x)] \]
from $X$ to $\Cu(A)$, along with, for each $[p] \in V(A)$, the value of 
\[ \langle a|_{\{x: [a(x)]=[p]\}} \rangle \]
in $V_c(\{x: [a(x)]=[p]\},A)$.
In fact, if we define
\[ V_c^{[p]}(Y,A) := \{\langle q \rangle \in V_c(Y,A): [q(x)]=[p]\ \forall x \in Y\} \]
then we see that $\langle a|_{\{x: [a(x)]=[p]\}} \rangle \in V_c^{[p]}(\{x: [a(x)]=[p]\},A)$.

In this section, we answer the natural question of the range of this invariant.
The answer is that, all possible values occur -- that is, every $\ll$-lower semicontinuous function $f:X \to \Cu(A)$ arises as $x \mapsto [a(x)]$ for some $a \in C_0(X,A \tens \K)_+$, and moreover, given $\langle a_{[p]} \rangle \in V_c^{[p]}(f^{-1}([p]),A)$ for each $[p] \in V(A)$, we can additionally satisfy
\[ \langle a|_{f^{-1}([p])} \rangle = \langle a_{[p]} \rangle\ \forall [p] \in V(A). \]
We now state this formally.

\begin{thm}\alabel{ASH-DataAttained}
Let $A$ be a simple, unital, $\jsZ$-stable ASH algebra and let $X$ be a second countable, locally compact Hausdorff space.
Let us be given a map $f:X \to \Cu(A)$ which is lower semicontinuous with respect to $\ll$ and, for each $[p] \in V(A)$, some $\langle a_{[p]} \rangle \in V_c(f^{-1}(\{[p]\}), A)$ such that $[a_{[p]}(x)] = [p]$ for all $x \in f^{-1}([p])$.
Then there exists $[a] \in \Cu(C_0(X,A))$ such that $[a(x)] = f(x)$ for all $x \in X$ and $\langle a|_{f^{-1}([p])} \rangle = \langle a_{[p]} \rangle$ for all $p \in V(A)$.
\end{thm}

To prove this, we first establish an important property of $\Cu(A)$: that it has Riesz interpolation.
The theorem will be proven using this fact and some other preliminary results.

\subsection{Riesz interpolation in $\Cu(A)$\alabel{CuInterpSection}}

We show here that $\Cu(A)$ has the Riesz interpolation property.
This is a well-known property for partially ordered sets, which we recall now.

\begin{defn}
A partially ordered set $(S,\leq)$ has the \textbf{Riesz interpolation property} if whenever elements $a_1,a_2,c_1,c_2 \in S$ satisfy
\[ \begin{array}{c} a_1 \\ a_2 \end{array} \leq \begin{array}{c} c_1 \\ c_2, \end{array} \]
there exists $b \in S$ satisfying
\[ \begin{array}{c} a_1 \\ a_2 \end{array} \leq b \leq \begin{array}{c} c_1 \\ c_2. \end{array} \]
Such an element $b$ is called an \textbf{interpolant}.
\end{defn}

Let us use $Lsc(K,[0,\infty])$ (resp.\ $Lsc(K,[0,\infty))$) to denote the set of all lower semicontinuous affine functions from a Choquet simplex $K$ to $[0,\infty]$ (resp.\ $[0,\infty)$).
With pointwise ordering, these are ordered semigroups.

In \cite{BrownPereraToms,ElliottRobertSantiago}, it was shown that when $A$ is a $\jsZ$-stable $C^*$-algebra which is simple and finite (and therefore stably finite), the Cuntz semigroup of $A$ can be identified with
\begin{equation} V(A) \dunion (Lsc(T(A),[0,\infty]) \setminus \{0\}). \alabel{ZstableCu} \end{equation}
The identification is made by identifying $[a] \in \Cu(A) \setminus V(A)$ with the function
\[ \tau \mapsto d_\tau(a) := \lim_n \tau(a^{1/n}). \]
The order is given as follows.
\begin{itemize}
\item For $[a] \in \Cu(A) \setminus V(A)$ and $[b] \in \Cu(A)$, we have $[a] \leq [b]$ if and only if $d_\tau(a) \leq d_\tau(b)$ for all $\tau \in T(A)$.
\item For $[p] \in V(A)$ and $[b] \in \Cu(A)$, we have $[p] < [b]$ if and only if $d_\tau(p) < d_\tau(b)$ for all $\tau \in T(A)$.
\end{itemize}
In particular, the induced order on $Lsc(T(A),[0,\infty])$ is the pointwise ordering.

Our proof that $\Cu(A)$ has Riesz interpolation will rest on the fact that $Lsc(T(A),\R)$ does (and thus, so does the convex subset $Lsc(T(A),[0,\infty))$).

\begin{lem}\alabel{Lsc-Interpolation}
Let $K$ be a Choquet simplex.
Then $Lsc(K,\R)$ has Riesz interpolation.
\end{lem}

\begin{proof}
Let us be given
\[ \begin{array}{c} a_1 \\ a_2 \end{array} \leq \begin{array}{c} c_1 \\ c_2 \end{array} \]
in $Lsc(K,\R))$.

By \ccite{Theorem 11.8}{Goodearl:book}, we have that $a_1$ is the pointwise supremum of all continuous affine functions $f:K \to \R$ which are below $a_1$ (and likewise of course for $a_2, c_1,$ and $c_2$).
For any finite set $f_1,\dots,f_k$ of continuous affine functions $f:K \to \R$ for which $f_i(x) \leq a_1(x)$ for all $x \in K$, consider the function $x \mapsto \max\{f_1(x),\dots,f_k(x)\}$.
A simple computation shows that this function is convex, so by \cite{Edwards}, there exists a continuous affine function $g:K \to \R$ such that
\[ \begin{array}{c} f_1 \\ \vdots \\ f_k \end{array} \leq g \leq a_1. \]
Thus, we may find an increasing net $(f_{\alpha,1})_{\alpha \in I_1}$ of continuous affine functions such that $a_1 = \sup f_{\alpha,1}$; in fact, the index set $I_1$ may be identified with the collection of finite subsets of continuous affine functions below $a_1$, and so, it has the property that for $\alpha \in I_1$, there are only finitely many $\beta \in I_1$ which come before $\alpha$.
Likewise, let $(f_{\alpha,2})_{\alpha \in I_2}$ be an increasing net of continuous affine functions such that $a_2 = \sup f_{\alpha,2}$, and let us assume that each $\alpha \in I_2$ has only finitely many indices before it.

We shall find an increasing net $(g_{\alpha_1,\alpha_2})_{I_1 \times I_2}$ of continuous affine functions, such that
\[ \begin{array}{c} f_{\alpha_1,1} \\ f_{\alpha_2,2} \end{array} \leq g_{\alpha_1,\alpha_2} \leq \begin{array}{c} c_1 \\ c_2 \end{array}. \]
Then, we shall define $b := \sup g_{\alpha_1,\alpha_2}$.
Since $b$ is the supremum of continuous functions, it is lower semicontinuous.
Moreover, since the net $(g_{\alpha_1,\alpha_2})$ is increasing, we see that $b$ is affine, so this interpolant is within $Lsc(K,\R)$, as required.

We define $g_{\alpha_1,\alpha_2}$ inductively, so that when defining $g_{\alpha_1,\alpha_2}$, we have already defined $g_{\beta_1,\beta_2}$ for all $(\beta_1,\beta_2) < (\alpha_1,\alpha_2)$.
We have that
\[ \begin{array}{c} f_{\alpha_1,1} \\ f_{\alpha_2,2} \\ g_{\beta_1,\beta_2} \end{array} \leq c_1, \]
for all $(\beta_1,\beta_2) < (\alpha_1,\alpha_2)$, so invoking \cite{Edwards} as above, we can find a continuous affine function $c'_1$ such that
\[ \begin{array}{c} f_{\alpha_1,1} \\ f_{\alpha_2,2} \\ g_{\beta_1,\beta_2} \end{array} \leq c'_1 \leq c_1, \]
for all $(\beta_1,\beta_2) < (\alpha_1,\alpha_2)$.
Likewise, we can find a continuous affine function $c'_2$ such that 
\[ \begin{array}{c} f_{\alpha_1,1} \\ f_{\alpha_2,2} \\ g_{\beta_1,\beta_2} \end{array} \leq c'_2 \leq c_2, \]
for all $\beta_1 < \alpha_1$ and all $\beta_2 < \alpha_2$.
By \cite{Edwards}, the group of continuous affine functions on $K$ has Riesz interpolation, and therefore, we may find a continouous affine function $g_{\alpha_1,\alpha_2}$ satisfying
\[ \begin{array}{c} f_{\alpha_1,1} \\ f_{\alpha_2,2} \\ g_{\beta_1,\beta_2} \end{array} \leq g_{\alpha_1,\alpha_2} \leq \begin{array}{c} c'_1 \\ c'_2, \end{array}, \]
for all $(\beta_1,\beta_2) < (\alpha_1,\alpha_2)$.
In particular, we have
\[ g_{\alpha_1,\alpha_2} \leq \begin{array}{c} c_1 \\ c_2, \end{array} \]
as required.
\end{proof}

\begin{prop}\alabel{Cu-Interpolation}
Let $A$ be a simple $\jsZ$-stable $C^*$-algebra.
Then $\Cu(A)$ has Riesz interpolation.
\end{prop}

\begin{proof}
Let us be given
\begin{equation}
\begin{array}{c} {[}a_1] \\ {[}a_2] \end{array} \leq \begin{array}{c} {[}c_1] \\ {[}c_2] \end{array}
\alabel{InterpRequired}
\end{equation}
in $\Cu(A)$.

It will suffice to find an interpolant under the assumption that $d_\tau(c_i) < \infty$ for all $\tau \in T(A)$; for then, in the general situation, we may approximate $[a_1],[a_2],[c_1],[c_2]$ by elements satisfying this restriction and \eqref{InterpRequired}, and obtain an increasing sequence of interpolants, the supremum of which is a bone fide interpolant for \eqref{InterpRequired}.
(The argument just outlined is essentially the same as in the proof of \ref{Lsc-Interpolation}, where interpolation in $Lsc(K,\R)$ was reduced to interpolation in the group of continuous affine functions from $K$ to $\R$.)

We may also assume that, in fact, we have $[a_i] < [c_j]$ for all $i,j$, since otherwise we would automatically have an interpolant.
Let $\hat{c}_j \in Lsc(T(A),[0,\infty))$ be the function $\hat{c}_j(\tau) = d_\tau(c_j)$.
Using \eqref{ZstableCu} to identifying $Lsc(T(A),[0,\infty])$ with a subset of $\Cu(A)$, we see that $\hat{c}_j \leq [c_j]$ and, since $[a_i] < [c_j]$, that $[a_i] \leq \hat{c}_j$.

We shall now also define functions $\hat{a}_i \in Lsc(T(A),[0,\infty))$ such that we have
\begin{equation} \begin{array}{c} \hat{a}_1 \\ \hat{a}_2 \end{array} \leq \begin{array}{c} \hat{c}_1 \\ \hat{c}_2 \end{array}, \alabel{InterpRequired2} \end{equation}
and $[a_i] \leq \hat{a}_i$.
To do this, when $[a_i] \not\in V(A)$, we set $\hat{a}_i = [a_i]$ (ie.\ the same definition as used for $\hat{c}_j$).
When $[a_i] \in V(A)$, we note that $[a_i] \ll \hat{c}_j = \sup_{\gamma < 1} \gamma\hat{c}_j$, and so there exists $\gamma \in (0,1)$ such that $[a_i] \leq \gamma\hat{c}_j$, for $j=1,2$.
Then, define $\hat{a}_i(\tau) = \gamma^{-1} d_\tau(a_i)$.

Now that we have \eqref{InterpRequired2}, it follows since $Lsc(T(A),\R)$ has Riesz interpolation that there exists $\hat{b} \in Lsc(T(A),\R)$ such that
\[ \begin{array}{c} \hat{a}_1 \\ \hat{a}_2 \end{array} \leq \hat{b} \leq \begin{array}{c} \hat{c}_1 \\ \hat{c}_2 \end{array},  \]
and thus, $\hat{b} = [b] \in \Cu(A)$ satisfies
\[ \begin{array}{c} {[}a_1] \\ {[}a_2] \end{array} \leq [b] \leq \begin{array}{c} {[}c_1] \\ {[}c_2] \end{array} \]
\end{proof}

\subsection{The proof of Theorem \ref{ASH-DataAttained}}

The proof of Theorem \ref{ASH-DataAttained} is done in a number of steps.
First, in Lemma \ref{RSH-DataAttained}, we show the existence of certain restricted elements in $\Cu(C_0(X,R))$ for recursive subhomogeneous $C^*$-algebras $R$; the restriction on the elements is in part related to the combined dimension of $X$ and of the total space of $R$ (both of which are assumed to be finite).
Using this result on the building blocks, Lemma \ref{ASH-RestrictedDataAttainment} establishes the result of Theorem \ref{ASH-DataAttained} in the special case that the range of $f$ is finite.
The proof of Theorem \ref{ASH-DataAttained} follows, combining Lemma \ref{ASH-RestrictedDataAttainment} and the Riesz interpolation property for $\Cu(A)$.

\begin{lem}\alabel{RSH-DataAttained}
Let $X$ be a finite dimensional locally compact Hausdorff space.
Let $R$ be a recursive subhomogeneous $C^*$-algebra with finite dimensional total space $\Omega$, and let $\sigma:R \to C(\Omega,\K)$ be the canonical representation of $R$.
Suppose that we are given 
\begin{enumerate}
\item[(i)] an open cover $U_1,\dots,U_n$ of $X$, such that each set $U_i$ is $\sigma$-compact.
\item[(ii)] for each $i=1,\dots,n$, an element $[a_i] \in \Cu(C(\overline{U_i},R))$.
\end{enumerate}
Suppose that, if $i \leq j$ then for $x \in \overline{U_i} \cap \overline{U_j}$ and $\omega \in \Omega$,
\begin{equation}
\rank \sigma(a_i(x))(\omega) + \frac{\dim X + d(\omega) - 1}2 \leq \rank \sigma(a_j(x))(\omega).
\alabel{RSH-DataAttained-RankDiffs}
\end{equation}
Then there exists $[a] \in \Cu(C(X,R))$ such that, for each $i$, if $f \in C_0(U_i \setminus \bigcup_{j > i} U_j)_+$ is strictly positive, then
\begin{equation}
[fa|_{U_i \setminus \bigcup_{j > i} U_j}] = [fa_i|_{U_i \setminus \bigcup_{j > i} U_j}]
\alabel{RSH-DataAttained-OpenCuntz}
\end{equation}
in $\Cu(C_0(U_i \setminus \bigcup_{j > i} U_j, R))$.
\end{lem}

\begin{rmk}
As seen in Proposition \ref{CuntzForC0Projs}, if $a_i$ is a projection then \eqref{RSH-DataAttained-OpenCuntz} amounts to
\[ \langle a|_{U_i \setminus \bigcup_{j > i} U_j} \rangle = \langle a_i|_{U_i \setminus \bigcup_{j > i} U_j} \rangle \]
in $V_c(U_i \setminus \bigcup_{j > i} U_j,R)$.
\end{rmk}

\begin{proof}
We shall find elements $s_i \in C(\overline{U_i},R)$ such that $s_i^*s_i = a_i$ and, for $i \leq j$,
\[ s_is_i^*|_{\overline{U_i} \cap \overline{U_j}} \in \her{s_js_j^*|_{\overline{U_i} \cap \overline{U_j}}}. \]
Then, using a strictly positive element $\lambda_i$ of $C_0(U_i)_+$ for each $i$, we shall set
\[ a = \sum_{i=1^n} \lambda_i s_is_i^*. \]
It follows that, if $f \in C_0(U_i \setminus \bigcup_{j > i} U_j)_+$ then
\[ fa|_{U_i \setminus \bigcup_{j > i} U_j} = f\lambda_i s_is_i^* + c \]
where $c \in \her{fs_is_i^*}$, and so
\[ [fa|_{U_i \setminus \bigcup_{j > i} U_j}] = [fs_is_i^*|_{U_i \setminus \bigcup_{j > i} U_j}] = [fa_i|_{U_i \setminus \bigcup_{j > i} U_j}]. \]

To find the elements $s_i$, we use induction on $i$, beginning at $n$ and decreasing.
For $i=n$, we simply set $s_n = a_n^{1/2}$.
Having defined $s_n,\dots,s_{i+1}$, let now define $s_i$.
This will be done first on $\overline{U_i} \cap \overline{U_{i+1}}$, and then extended to add the set $\overline{U_i} \cap \overline{U_{i+2}}$, and so on until we add the set $\overline{U_i} \cap \overline{U_n}$, and then finally to the rest of $\overline{U_i}$.

For the step where we extend the definition to include the set $\overline{U}_i \cap \overline{U}_j$ (for $j > i$), we can assume that $s_i$ is already defined on some closed (possibly empty) subset $K$ of $\overline{U_i} \cap \overline{U_j}$, such that the definition already satisfies $s_i \in \her{s_j}$.
By \eqref{RSH-DataAttained-RankDiffs}, we can apply Lemma \ref{RSH-Embedding}, providing the extension of $s_i$ to $\overline{U_i} \cap \overline{U_j}$, as required.
\end{proof}

The following was used, implicitly, in the proof of \ccite{Proposition 4.6}{CommutativeCuntz}, in the case that $A=\K$.

\begin{lem} \alabel{ProjNbhdExtension}
Let $X$ be a locally compact Hausdorff space and let $Y$ be a closed subset.
Let $A$ be a $C^*$-algebra.
Let $p \in C_b(Y,A)$ such that $p(x)$ is a projection for all $x \in Y$, and $p(x) \mvneq p(y)$ for all $x,y \in Y$. 
Then there exists an open set $U$ with $Y \subseteq U \subseteq X$ and some $\tilde{p} \in C_b(U,A)$ such that $\tilde{p}|_Y = p$, $\tilde{p}(x)$ is a projection for all $x \in U$ and $\tilde{p}(x) \mvneq \tilde{p}(y)$ for all $x,y \in U$.
\end{lem}

\begin{proof}
We may find a continuous extension $a \in C_b(X,A)$ of $p$.
Moreover, we may find an open set $U$ containing $Y$ such that, for every $x \in \overline{U}$ there exists $y \in Y$ such that $\|a(x)-a(y)\|<1/2$.
It follows that the spectrum of $a|_U$ (in the algebra $C_b(U,A)$) is contained in $\R \setminus \{1/2\}$, and so, if
\[ f(t) = \begin{cases} 0,\ &\text{for }t<1/2, \\ 1,\ &\text{for }t>1/2, \end{cases} \]
then by functional calculus, $\tilde{p} := f(a|_U) \in C_b(U,A)$.
Also, clearly $\tilde{p}|_Y = p$ and, for $x \in U$ there exists $y \in Y$ such that $\tilde{p}(x) \mvneq p(y)$.
Consequently, the Murray-von Neumann class of $\tilde{p(x)}$ is constant over all of $U$.
\end{proof}

\begin{lem}\alabel{ASH-RestrictedDataAttainment}
Let $A$ be a simple, unital $\jsZ$-stable ASH algebra and let $X$ be a second countable, locally compact Hausdorff space.
Let us be given a map $f:X \to \Cu(A)$ which is lower semicontinuous with respect to $\ll$ and, for each $[p] \in V(A)$, some of $\langle a_{[p]} \rangle \in V_c(f^{-1}(\{[p]\}), A)$ such that $[a_{[p]}(x)] = [p]$ for all $x \in f^{-1}([p])$.
If the range of $f$ is finite then there exists $[a] \in \Cu(C_0(X,A))$ such that $[a(x)] = f(x)$ for all $x \in X$ and $\langle a|_{f^{-1}([p])} \rangle = \langle a_{[p]} \rangle$ for all $p \in V(A)$.
\end{lem}

\begin{proof}
Let $[b_1],\dots,[b_n] \in \Cu(A)$ be the values taken by $f$, in non-decreasing order.
Let us first prove something weaker.
Namely, if we in fact have $\e > 0$ such that, whenever $[b_i] < [b_j]$,
\[ [b_i] \leq [(b_j-\e)_+], \]
then we can find $[a] \in \Cu(C_0(X,A))$ such that, for each $x$, if$f(x) = [b_i]$, then
\[ [(b_i-\e)_+] \leq [a(x)] \leq [b_i], \]
and for each $i$ for which $[b_i] \in V(A)$, we still have
\[ \langle a|_{f^{-1}([b_i])} \rangle = \langle a_{[b_i]} \rangle. \]

To prove this, first let us define open sets $U_i$ as follows.
For each $i$ for which $[b_i] \not\in V(A)$, set $U_i := \{x: f(x) \geq [b_i]\}$; by (i) and by $\ll$-lower semicontinuity of $f$, we see that $U_i$ is open.
For each $i$ for which $[b_i] \in V(A)$, use Lemma \ref{ProjNbhdExtension} to find $U_i \subseteq \{x: f(x) \geq [b_i]\}$ upon which $a_{[b_i]}$ extends continuously (and such that $[a_{[b_i]}(x)] = [b_i]$ for all $x \in U_i$).

We may write $X$ as an inverse limit of finite-dimensional locally compact Hausdorff spaces
\[ X_1 \labelledleftarrow{\alpha_2^1} X_2 \labelledleftarrow{\alpha_3^2} \cdots \leftarrow X, \]
where the connecting maps, $\alpha_j^i$, are proper and surjective.
For each $i$, we may find open sets $U_{i,k} \subseteq X_k$ such that $(\alpha_k^{k+1})^{-1}(U_{i,k}) \subseteq U_{i,k+1}$ and
\[ U_i = \bigcup_k (\alpha_k^{\infty})^{-1}(U_{i,k}). \]
In fact, we may assume that whenever $[b_i] \in V(A)$, the set $(\alpha_k^\infty)^{-1}(U_{i,k})$ is compactly contained in $U_k$, and (by possibly replacing $(X_k)$ with a subsequence) that there exists $[p_{i,k}] \in V(C(\overline{U_{i,k}},A))$ such that
\[ [p_{i,k} \circ \alpha_k^\infty|_{(\alpha_k^\infty)^{-1}(\overline{U_{i,k}})}] = [p_i|_{(\alpha_k^\infty)^{-1}(\overline{U_{i,k}})}]. \] 

For each $k$, we shall find $[c_k] \in \Cu(C_0(X_k,A))$ such that, for each $x \in U_{i,k} \setminus \bigcup_{j > i} U_{j,k}$, we have
\[ [(b_i-\e)_+] \leq [c_k(x)] \leq [b_i], \]
and whenever $[b_i] \in V(A)$, we have
\[ \langle c_k|_{U_{i,k} \setminus \bigcup_{j > i} U_{j,k}} \rangle = \langle p_{i,k}|_{ U_{i,k} \setminus \bigcup_{j > i} U_{j,k}} \rangle. \]
It follows, using Theorem \ref{ASH-Comparison}, that $[c_k \circ \alpha_k^{k+1}] \leq [c_{k+1}]$, and it is easy to see that 
\[ [a] := \sup [c_k \circ \alpha_k^\infty] \]
is as needed.

We describe now how to find $[c_k]$.
Let
\[ A_1 \labelledrightarrow{\phi_1^2} A_2 \labelledrightarrow{\phi_2^3} \cdots \rightarrow A \]
be a directed system as in Proposition \ref{ASH-SlowCharacterization} (iv).
We may therefore find some $n$ such that:
\begin{enumerate}
\item[(i)] Within $\Cu(A_n)$, there are elements $[\hat{b}_i]$ such that
\[ [(b_i-\e/2)_+] \leq [\phi_n^\infty(\hat{b}_i)] \leq [b_i], \text{ and} \]
\item[(ii)] Within $V(C(\overline{U_{i,k}},A_n))$ there is an element $[\hat{p}_{i,k}]$ such that
\[ [(\id_{C(\overline{U_{i,k}})} \tens \phi_n^\infty)(\hat{p}_{i,k})] = [p_{i,k}]. \]
(For such $i$, we may assume that $[\hat{b}_i] = [\hat{p}_{i,k}(x)]$, which is the same for all $x$.)
\end{enumerate}

By (i) and \ccite{Proposition 2.2}{Rordam:UHFII}, let $\dl > 0$ be such that $[(b_i-\e)_+] \leq [(\phi_n^\infty(\hat{b}_i)-\dl)_+]$ for all $i$.
We can see that, by possibly increasing $n$, whenever $[b_i] < [b_j]$, it is the case that
\[ [\hat{b}_i] < [(\hat{b}_j-\dl)_+]. \]
But then, for such $i,j$, using the same argument as in the claim in the proof of Lemma \ref{ConditionToBuildingBlocks} (iii), there exists a non-zero $c$ such that
\[ (\hat{b}_i-\dl)_+ \dsum c \culeq (\hat{b}_j-\dl)_+ \]
and so, using Proposition \ref{ASH-SlowCharacterization} (iv) as in the proof of Lemma \ref{ConditionToBuildingBlocks} (iii), by increasing $n$, we may assume that for all $\omega$ in the total space of $A_n$,
\[ \rank \ev_\omega (\hat{b}_i-\dl)_+ + \frac{\dim(X_k) + d(\omega) - 1}2 \leq \rank \ev_\omega (\hat{b}_j-\dl)_+. \]
We see that Lemma \ref{RSH-DataAttained} applies, providing us with $[\hat{c}_k] \in \Cu(C_0(X_k,A_n))$ such that, for $x \in U_{i,k} \setminus \bigcup_{j > i} U_{j,k}$,
\[ [\hat{c}_k(x)] = [(\hat{b}_i-\dl)_+], \]
and when $[b_i] \in V(A)$,
\[ \langle \hat{c}_k|_{U_{i,k} \setminus \bigcup_{j > i} U_{j,k}} \rangle = \langle \hat{p}_{i,k}|_{U_{i,k} \setminus \bigcup_{j > i} U_{j,k}} \rangle. \]
Hence, $[c_k] = [(\id_{C_0(X_k)} \tens \phi_n^\infty)(\hat{c}_k)]$ is as required.

Having completed the proof of the weaker statement, let us now prove the full theorem.
Let $\dl_{1,1} > 0$.
Using \ccite{Proposition 2.2}{Rordam:UHFII}, we may iteratively find $\dl_{i,1} > 0$ such that, whenever $[b_i] \leq [b_j]$, we have
\[ [(b_i-\dl_{i,1}/2)_+] \leq [(b_j-\dl_{j,1})_+]. \]
Whenever $[b_i] \in V(A)$, we may (by picking $\dl_i$ small enough) assume that $[(b_i-\dl_i)_+] = [b_i]$.
Then, by applying the restricted result above, we can find $[c_1] \in \Cu(C_0(X,A))$ such that whenever $f(x) = [b_i]$, we have
\[ [(b_i-\dl_{i,1})_+] \leq [c_1(x)] \leq [(b_i-\dl_{i,1}/2)_+], \] 
and $\langle c_1|_{f^{-1}([p])} \rangle = \langle a_{[p]} \rangle$ for all $[p] \in V(A)$.
(Note that, strictly speaking, we may not quite satisfy the hypothesis of the restricted result above, since we may have that $[(b_i-\dl_{i,1}/2)_+] \leq [(b_j-\dl_{j,1}/2)_+]$ even though $[b_i] \not\leq [b_j]$.
But, looking at how the hypothesis is used, we see that the argument above is not affecting by this potential shortcoming.)

Next, we let $\dl_{1,2} \in (0,\dl_{1,1}/2)$ and, working inductively, select $\dl_{i,2} \in (0,\dl_{i,1}/2)$ such that, whenever $[b_i] \leq [b_j]$,
\[ [(b_i-\dl_{i,2}/2)_+] \leq [(b_j-\dl_{j,2})_+]. \]
Thus, we obtain $[c_2] \in \Cu(C_0(X,A))$ such that whenever $f(x) = [b_i]$, we have
\[ [(b_i-\dl_{i,2})_+] \leq [c_1(x)] \leq [(b_i-\dl_{i,2}/2)_+], \]
and $\langle c_1|_{f^{-1}([p])} \rangle = \langle a_{[p]} \rangle$ for all $[p] \in V(A)$.
Since $\dl_{i,2} < \dl_{i,1}/2$, Theorem \ref{ASH-Comparison} ensures that $[c_1] \leq [c_2]$.
By continuing, we form an increasing sequence $([c_k]) \subset \Cu(C_0(X,A))$ such that, for each $x \in X$, $f(x) = \sup [c_k(x)]$, and, for all $[p] \in V(A)$,
\[ \langle c_k|_{f^{-1}([p])} \rangle = \langle a_{[p]} \rangle. \]

It follows that $[a] := \sup [c_k]$ is an element satisfying the conclusion of this lemma.
\end{proof}

\begin{proof}[Proof of Theorem \ref{ASH-DataAttained}]
For $[p] \in V(A)$, using Lemma \ref{ProjNbhdExtension}, let $U_{[p]}$ be open such that
\[ f^{-1}([p]) \subseteq U_{[p]} \subseteq \{x: f(x) \geq [p]\}, \]
such that $a_{[p]}$ extends continuously, as a projection, to $U_{[p]}$ (and for $x \in U_{[p]}$, we have $[a_{[p]}(x)] = [p]$).
For $[b] \not\in V(A)$, set
\[ U_{[b]} := \{x \in X: [b] \ll f(x)\}, \]
which is an open set.

We will want to approximate $f$ by an increasing sequence of functions $f_k$ which have finite range, then use Lemma \ref{ASH-RestrictedDataAttainment} with $f_k$.
Let $(S_k)$ be an increasing sequence of finite subsets of $\Cu(A)$ such that every $[b] \in \Cu(A)$ is the supremum of all $[t] \in \bigcup_k S_k$ for which $[t] \ll [b]$.
Ideally, we could set
\[ f_k(x) = \sup \{[t] \in S_k: x \in U_{[t]}\}, \]
which would work if $\Cu(A)$ was lattice-ordered.
Although $\Cu(A)$ may not be lattice-ordered, we do know that $\Cu(A)$ has the Riesz interpolation property, so instead we shall want to set $f_k(x)$ to be an interpolant between the elements $[t]$ and $f(x)$; this is the object of the following argument.

For each $[b] \in \Cu(A)$, we may find an increasing sequence $(U_{[b],k})_{k=1}^\infty$ of open sets of $U_{[b]}$, such that each $U_{[b],k}$ is compactly contained in $U_{[b]}$ and
\[ U_{[b]} = \bigcup_k U_{[b],k}. \]

Let us first construct $f_1$.
For every subset $T$ of $S_1$, let $\max(T)$ denote the maximal elements of $T$ (that is, all $[t] \in T$ such that there is no $[s] \in T$ for which $[t] < [s]$).
If $T=\max(T)$, then this means that the elements of $T$ are pairwise incomparable.
In this case, we set
\[ W_{T,1} = \bigcap_{[t] \in T} U_{[t],1}. \]
We shall associate an element $[b_T] \in \Cu(A)$ to such $T$, and these elements shall satisfy
\begin{enumerate}
\item[(i)] $[b_\emptyset] = [0]$ and $[b_T] = [t]$ if $T = \{[t]\}$;
\item[(ii)] If $T$ is not a singleton then $[b_T] \not\in V(A)$;
\item[(iii)] If $T_1,T_2$ are both subsets of $S_1$ and every element of $T_1$ is dominated by some element of $T_1$ then
\[ [b_{\max(T_1)}] \leq [b_{\max(T_2)}]; \text{ and} \]
\item[(iv)] For each subset $T$, if $x \in W_{\max(T),1}$, then $[b_T] \ll f(x)$.
\end{enumerate}

First, set $[b_{[t]}] = [t]$ to satisfy (i); by the choice of $U_{[t],1}$, we can see that (iii) holds for $T=\{[t]\}$.
We shall iteratively choose the remaining elements $[b_T]$, in order that we satisfy (ii)-(iv).
Note that, when $T_1,T_2$ are subsets of $S_1$, we have that every element of $T_1$ is dominated by some element of $T_2$ if and only if the same relation holds for $\max(T_1)$ and $\max(T_2)$.
Let us denote this relation by $\preccurlyeq$; we shall iterate through the subsets $T$ of $S_1$ which have pairwise incomparable elements, in non-decreasing order under $\preccurlyeq$.
Setting $W'_T$ to be the union of every set $W_{T_1}$ for which $T \preccurlyeq T_1$, we shall strengthen (iv) to:
\begin{enumerate}
\item[(iv)$'$] $[b_T] \ll f(x)$ for all $x \in \overline{W'_T}$.
\end{enumerate}

Let us see now how to get the element $[b_T]$.
By the choice of the sets $U_{[t],1}$, we have that $\overline{W'_T}$ is compact.
Let $x \in \overline{W'_T}$.
By induction, we have that $[b_{T_0}] \ll f(x)$ for all $T_0 \preccurlyeq T$.
Since $f(x)$ is a supremum of an increasing sequence of elements that are $\ll f(x)$, there exists $[s_x] \in \Cu(A)$ such that
\[ [b_{T_0}] \leq [s_x] \ll f(x), \]
for all $T_0 \preccurlyeq T$.
Since $f$ is lower semicontinuous with respect to $\ll$, there is an open neighbourhood $V_x$ of $x$ such that $[s_x] \ll [y]$ for all $[y] \in V_x$.

Using compactness of $\overline{W'_T}$, there exists $[s_1],\dots,[s_\ell] \in \Cu(A)$ such that $[b_{T_0}] \leq [s_i]$ for all $T_0 \preccurlyeq T$ and all $i=1,\dots, \ell$, and for each $x \in W'_T$, $[s_i] \ll f(x)$ for some $i$.
By Proposition \ref{Cu-Interpolation}, $\Cu(A)$ has the Riesz interpolation, so we may find $[b_T]$ such that
\[ [b_{T_0}] \leq [b_T] \leq [s_i] \]
for all $T_0 \preccurlyeq T$ and all $i=1,\dots,\ell$.
Hence, we have $[b_T] \ll f(x)$ for all $x \in \overline{W'_T}$.
It is not hard to see that, if $[b_T] \in V(A)$ then we may replace $[b_T]$ by the element $[b] \in \Cu(A) \setminus V(A)$ for which $d_\tau(b) = d_\tau(b_T)$ for all $\tau \in T(A)$; thus, we can always obtain $[b_T] \not\in V(A)$.

Having now defined $[b_T]$ for all subsets $T$ of $S_1$ consisting of pairwise incomparable elements, we define $f_1$ as follows.
For $x \in X$, set $T(x) := \{[t] \in S_1: x \in U_{[t],1}\}$, and then set
\[ f_1(x) := [b_{T(x)'}]. \]
Equivalently (using condition (iii) above), we have that
\[ f_1(x) = \max \{[b_{\max(T)}]: x \in W_{\max(T),1}\}, \]
which shows that $f_1$ is lower semicontinuous.

The definition of $f_2$ is similar, but an additional measure of care is needed to ensure that $f_1(x) \leq f_2(x)$ for all $x$.
In order to arrange this, we set $S_2' = S_2 \cup \{[b_{\max(T)}]: T \subseteq S_1\}$ and, for $[t] \in S_2'$, we define $U'_{[t],2}$ as follows.
If $[t] = [b_{\max(T)}]$ for $T \subseteq S_1$ then we set
\[ U'_{[t],2} := U_{[t],2} \cup W_{\max(T),1}, \]
otherwise, set
\[ U'_{[t],2} := U_{[t],2}. \]
Now, we do the same construction as for $f_1$, except using $S_2'$ in place of $S_1$ and $U'_{[t],2}$ in place of $U_{[t],1}$.
This produces $f_2$, and continuing, we create a sequence $(f_k)$ of lower semicontinuous functions from $X$ to $\Cu(A)$, such that for each $x \in X$, $f(x)$ is the supremum of the increasing sequence $f_k(x)$.

Now, using Lemma \ref{ASH-RestrictedDataAttainment}, we may find, for each $k$, some $[b_k] \in \Cu(C_0(X,A))$ such that $[b_k(x)] = f_k(x)$ for all $x$ and, for each $[p] \in V(A)$,
\[ \langle b_k|_{f_k^{-1}([p])} \rangle = \langle a_{[p]}|_{f_k^{-1}([p])} \rangle \]
(remember that $a_{[p]}$ extends continuously to $U_{[p]}$, and our construction of $f_k$ ensures that $f_k^{-1}([p]) \subseteq U_{[p]}$).
Hence, by Theorem \ref{ASH-Comparison}, the sequence $[b_k]$ is increasing, and evidently its supremum, $[a]$, satisfies the conclusions of this lemma.
\end{proof}
\section{The equivalence of $Ell(A)$ and $\Cu(C(\mathbb{T},A))$ \alabel{Ell-CuT-Eq}}

An application of the Cuntz semigroup computation in Theorem \ref{MainThm} arises by considering the case that $X$ is the circle; we see that, for the algebras in Theorem \ref{MainThm}, the Cuntz semigroup of $C(\mathbb{T},A)$ contains the same information as the Elliott invariant.
Recall that the Elliott invariant of a simple, unital $C^*$-algebra $A$ is
\[ Ell(A) := (K_0(A), K_0(A)_+, [1]_{K_0(A)}, K_1(A), T(A), \rho_A). \]
Here, $K_0(A)$ and $K_1(A)$ are the $K$-groups of the $C^*$-algebra $A$, while $K_0(A)_+$ is the positive cone of $K_0(A)$ (this is the subsemigroup of $K_0(A)$ generated by the images of projections in $A \tens \K$) and $[1]_{K_0(A)}$ is the image of $1_A$ in $K_0(A)$.
As before, we use $T(A)$ to denote the Choquet simplex of tracial states on $A$ and $\rho_A:K_0(A) \times T(A) \to \R$ denotes the pairing given by
\[ \rho_A([p],\tau) := \tau(p). \]

It is known, largely as a consequence of the description \eqref{ZstableCu} of $\Cu(A)$, that for simple, finite, exact $\jsZ$-stable $C^*$-algebras, $(\Cu(A),[1]_{\Cu(A)})$ is functorially equivalent to $(K_0(A),K_0(A)_+,[1]_{K_0(A)},T(A),\rho_A)$.
Let us review the argument.
Given $(\Cu(A),[1]_{\Cu(A)})$, we obtain $V(A)$ as the subsemigroup of elements $[a] \in \Cu(A)$ for which $[a] \ll [a]$ (Proposition \ref{BCmagic}).
From this, we construct $(K_0(A),K_0(A)_+)$.
It was shown in \cite{BlackadarHandelman} that the map $\tau \to d_\tau$ is a bijection between the normalized $2$-quasitraces on $A$ and the lower semicontinuous, additive, order-preserving functionals on $\Cu(A)$ which map $[1]$ to $1$.
Combining this with Uffe Haagerup's unpublished result that $2$-quasitraces are traces on a unital exact $C^*$-algebra \cite{Haagerup:quasitraces} (see also \cite{BrownWinter:quasitraces} for a less general result, which nonetheless applies to exact $\jsZ$-stable algebras) yields that $T(A)$ can be obtained from $(\Cu(A),[1]_{\Cu(A)})$; the pairing $\rho_A$ arises naturally also.
The opposite direction follows immediately from the description \eqref{ZstableCu} together with the fact that $V(A) \iso K_0(A)_+$; this fact amounts to the statement that the projections in $A \tens \K$ satisfy cancellation, and to see this, we find that $A$ has stable rank one by \cite{Rordam:Z} and thus cancellation follows by \ccite{Proposition 6.5.1}{Blackadar:KBook}.

Evidently, if we want to recover the Elliott invariant, what is missing in the Cuntz semigroup is the $K_1$-group.
However, $\Cu(C(\mathbb{T},A))$ contains the Murray-von Neumann semigroup of $C(\mathbb{T},A)$, and we show (Proposition \ref{ProjCancellation} and paragraphs following) that for simple, unital, $\jsZ$-stable ASH algebras, this is simply
\[ K_0(C(\mathbb{T},A))_+ \iso \{0\} \dunion K_0(A)_{++} \times K_1(A). \]
Here, and in what follows, $K_0(A)_{++}$ denotes the strictly positive cone, ie.\ it is all of $K_0(A)_+$ except $0$.
From this, we can recover $K_1(A)$ as the subgroup of $K_0(C(\mathbb{T},A))$ consisting of elements $g$ such that
\[ g + h \geq 0 \text{ and } g-h \leq 0 \text{ whenever }h \in K_0(C(\mathbb{T},A))_{++}. \]
Of course, we can also obtain $\Cu(A)$ from $\Cu(C(\mathbb{T},A))$, since it is isomorphic to the quotient of $\Cu(C(\mathbb{T},A))$ by any maximal ideal $I$ (where by an ideal we mean a hereditary subsemigroup closed under increasing sequential suprema; see \cite{CiupercaRobertSantiago}).

The application of Theorem \ref{MainThm} is to show that we can go the other way: given $Ell(A)$, we can recover $\Cu(C(\mathbb{T},A))$.
Indeed, it is shown (in Proposition \ref{CuT-Descr}) that, when $X=\mathbb{T}$ in Theorem \ref{MainThm}, then the only time that the elements $\langle q_{[p]} \rangle \in V_c(f^{-1}([p]),A)$ can give non-trivial information is when $f^{-1}([p]) = \mathbb{T}$, ie.\ $f$ takes a constant value in $V(A)$.
Combining with the above analysis of $V(C(\mathbb{T},A))$, this yields
\begin{align}
\alabel{CuT-DescrEq}
\Cu(C(\mathbb{T},A)) &\iso \{f:\mathbb{T} \to K_0(A)_{++} \dunion Lsc(T(A),[0,\infty]): \\
\notag &\qquad f \text{ is } \ll\text{-lower semicontinuous, } f \not\equiv [p] \in K_0(A)_{++}\} \\
\notag & \quad \dunion (K_0(A)_{++} \times K_1(A)).
\end{align}

The equivalence of the invariants $Ell(A)$ and
\[ \Cu_{\mathbb{T}}(A) := (\Cu(C(\mathbb{T}, A)), [1]_{\Cu(C(\mathbb{T},A))}) \]
established is functorial, in the following sense.
Given a unital $*$-homomorphism $\phi:A \to B$, this induces morphisms $\phi_1 := [\phi]_{Ell(A)}$ and $\phi_2 := [\phi]_{\Cu_{\mathbb{T}}(A)}$ of the two invariants, and we have the following.

\begin{prop}
The map $\phi_1$ induces a morphism  $(\phi_1)_*:\Cu_{\mathbb{T}}(A) \to \Cu_{\mathbb{T}}(B)$ and $\phi_2$ induces a morphism $(\phi_2)_*:Ell(A) \to Ell(B)$, and we have $(\phi_1)_* = \phi_2$ and $(\phi_2)_* = \phi_1$.
\end{prop}

\begin{proof}
The map $\phi_1:Ell(A) \to Ell(B)$ consists of maps $\phi_1^{K_0}:K_0(A) \to K_0(B)$, $\phi_1^{K_1}:K_1(A) \to K_1(B)$, and $\phi_1^T:T(B) \to T(A)$.
Moreover, $\phi_1^{K_0}$ sends $K_0(A)_+$ into $K_0(B)_+$ and sends $[1_A]$ to $[1_B]$.
Let us use \eqref{CuT-DescrEq} to describe the induced map $\Cu(C(\mathbb{T},A)) \to \Cu(C(\mathbb{T},B))$.

First, $\phi_1$ induces a map $\Cu(A) \to \Cu(B)$ by sending $K_0(A)_+$ to $K_0(B)_+$ and sending $f \in Lsc(T(A),[0,\infty])$ to $f \circ \phi_1^T$.
Now, given a $\ll$-lower semicontinuous map $f:\mathbb{T} \to \Cu(A)$, we compose with the map just described to get a $\ll$-lower semicontinuous map $(\phi_1)_*(f):\mathbb{T} \to \Cu(B)$.

\begin{claim} If the $\ll$-lower semicontinuous map $f:\mathbb{T} \to \Cu(A)$ does not take a constant value in $K_0(A)_{++}$ then $(\phi_1)_*(f)$ does not take a constant value in $K_0(B)_{++}$.
\end{claim}

\begin{proof}[Proof of claim.]
If $(\phi_1)_*(f) \equiv [p] \in K_0(B)_{++}$ then, first it is clear by the construction of $(\phi_1)_*(f)$ that $f(t) \in K_0(A)_{++}$ for all $t \in \mathbb{T}$.
Since $f$ is $\ll$-lower semicontinuous and not constant, it follows that there are projections $p,q \in A \tens \K$ such that $p < q$ but $\phi(p) \mvneq \phi(q)$.
Since $\phi$ is unital and its domain, $A$, is simple, it must be injective.
Thus we have $\phi(p) < \phi(q)$, yet $\phi(p) \mvneq \phi(q)$, which contradicts the fact that $B$ is stably finite.
\end{proof}

Given $(g_0,g_1) \in K_0(A)_{++} \times K_1(A)$, we set $(\phi_1)_*(g_0,g_1) = (\phi_1^{K_0}(g_0), \phi_1^{K_1}(g_1))$, if $\phi_1^{K_0}(g_0) \neq 0$; otherwise, set $(\phi_1)_*(g_0,g_1) = 0$.

This completes the description of $(\phi_1)_*$.
It is easily checked that $(\phi_1)_* = \phi_2$.

Going the other way, let us be given the map $\phi_2:\Cu(C(\mathbb{T},A)) \to \Cu(C(\mathbb{T},B))$.
If $I$ is a maximal ideal of $\Cu(C(\mathbb{T},A))$ then $\phi_2(I)$ generates some maximal ideal $J$ of $\Cu(C(\mathbb{T},B))$, and so $\phi_2$ induces a function
\[ \Cu(A) \iso \Cu(C(\mathbb{T},A))/I \to \Cu(C(\mathbb{T},B))/J \iso \Cu(B). \]
One can easily see that this does not depend on the choice of the maximal ideal $I$.
From this function we obtain a morphism
\[ (K_0(A),K_0(A)_+,[1]_{K_0(A)},T(A),\rho_A) \to (K_0(B),K_0(B)_+,[1]_{K_0(B)},T(B),\rho_B) \]
that agrees with $\phi_1$.

Lastly, $\phi_2$ restricts to a map from $K_0(C(\mathbb{T},A))_+$ to $K_0(C(\mathbb{T},B))_+$, and this induces a map
\[ \psi: K_0(A) \dsum K_1(A) \iso K_0(C(\mathbb{T},A)) \to K_0(C(\mathbb{T},B)) \iso K_0(B) \dsum K_1(B). \]

\begin{claim}
The map $\psi$ takes $K_1(A)$ into $K_1(B)$.
\end{claim}

\begin{proof}[Proof of claim.]
Let $g \in K_1(A)$; so $g$ can be represented as a difference $(g_0,g_1)-(g_0,g'_1)$ where $g_0 \in K_0(A)_{++}$ and $g_1,g'_1 \in K_1(A)$.
Then $\psi(g) = \phi_2(g_0,g'_1) - \phi_2(g_0,g'_1)$.
Let $I,J$ be maximal ideals of $\Cu(C(\mathbb{T},A)), \Cu(C(\mathbb{T},B))$ respectively as above, and $\pi_I:\Cu(C(\mathbb{T},A)) \to \Cu(C(\mathbb{T},A))/I, \pi_J: \Cu(C(\mathbb{T},B)) \to \Cu(C(\mathbb{T},B))/J$ are the quotient maps.
Note that $\pi_I(g_0,g_1) = g_0 = \pi_I(g_0,g'_1)$, and so
\[ \pi_J(\phi_2(g_0,g_1)) = (\phi_2)_*(g_0) = \pi_J(\phi_2(g_0,g'_1)). \]
It is easy to see that this ensures that $\phi_2(g_0,g_1)-\phi_2(g_0,g'_1) \in K_1(B)$.
\end{proof}

It is clear that the $K_1$-component of $(\phi_2)_*$ just described is the same as the $K_1$-component of $\phi_1$, so altogether $(\phi_2)_* = \phi_1$.
\sqbox 
\end{proof}

Functorial equivalence of $C^*$-algebra isomorphism invariants is relevant to the Elliott programme of classification of $C^*$-algebras.
We say that a class $\mathcal{C}$ of $C^*$-algebras is classified by the invariant $I$ (a functor from the category whose objects are those $C^*$-algebras in $\mathcal{C}$ and whose morphisms are $*$-homomorphisms) if, for any $C^*$-algebras $A,B \in \mathcal{C}$ and any isomorphism $\alpha:I(A) \to I(B)$ of the invariants, there exists an isomorphism $\phi:A \to B$ which lifts $\alpha$, ie.\ such that $I(\phi) = \alpha$.
The Elliott invariant, $Ell(\cdot)$, has been shown to classify substantial classes of simple, nuclear, unital, separable $C^*$-algebras, including $\jsZ$-stable approximately homogeneous $C^*$-algebras \cite{ElliottGongLi:AHclassification} (see also \ccite{Corollary 6.7}{Winter:Perfect}, where other descriptions of this class are given), and $\jsZ$-stable ASH algebras with slow dimension growth for which projections separate traces, \ccite{Corollary 5.5}{LinNiu:Lifting}.

As a consequence of the functorial equivalence of $Ell(\cdot)$ and $\Cu_{\mathbb{T}}(\cdot)$ for simple, unital, $\jsZ$-stable ASH algebras, we see that any subclass $\mathcal{C}$ of such algebras is classified by $Ell(\cdot)$ if and only if it is classified by $\Cu_{\mathbb{T}}(\cdot)$.
A comparable result is obtained in \cite{PereraToms:recasting,BrownPereraToms}; there, it is established that for simple, unital, exact $\jsZ$-stable algebras, $Ell(\cdot)$ is functorially equivalent to the invariant $A \mapsto (Ell(A),W(A), [1]_{W(A)})$.
A key difference is that, in that result, the Elliott invariant is augmented, rather than replaced, by a Cuntz semigroup invariant.

We now prove the technical points mentioned above.
First, we establish that projections in $C(\mathbb{T},A \tens \K)$ satisfy cancellation; in fact, we show this more generally, where $\mathbb{T}$ can be replaced by any metrizable compact Hausdorff space.

\begin{prop}\alabel{ProjCancellation}
Let $K$ be a metrizable compact Hausdorff space $K$ and let $A$ be a simple, unital, $\jsZ$-stable ASH algebr.
Then the projections in $C(K, A \tens \K)$ satisfy cancellation.
\end{prop}

\begin{proof}
Let $[e_1],[e_2],[f] \in V(C(K,A \tens \K))$ be elements satisfying
\[ [e_1] + [f] = [e_2] + [f]. \]
Since $[e_1(x)]$ are locally constant, we may decompose $K$ into finitely many components such that these are constant on each of them.
It will then suffice to show that $[e_1] = [e_2]$ holds for the restriction to each of these components.
This shows that we reduce the problem to the case where we assume that $[e_1(x)]$ takes a constant value.
If $[e_1] = 0$ then there is nothing to show, so let us assume otherwise.

Since we can write $K$ as an inverse limit of finite-dimensional spaces, let $K'$ be a finite-dimensional space, $\alpha:K \to K'$ an injective map, $[e'_1],[e'_2],[f'] \in V(C(K',A \tens \K))$ lifts of $[e_1],[e_2],[f]$ respectively, such that
\[ [e'_1] + [f'] = [e'_2] + [f']. \]
Using Proposition \ref{ASH-SlowCharacterization} (iv), we may find a recursive subhomogeneous algebra $R$ and lifts $[e''_1],[e''_2],[f''] \in V(C(K', R \tens \K))$ of $[e'_1],[e'_2],[f']$ such that
\[ [e''_1] + [f''] = [e''_2] + [f''], \]
and $\rank \ev_\omega(e''_1(x)) \geq (\dim(K') + d(\omega))/2$ for all $x \in K'$ and all $\omega$ in the total space of $R$.
It follows by \ccite{Theorem 4.6}{Phillips:rsh}, that $[e''_1] = [e''_2]$.
Returning to $C(K,A \tens \K)$ gives $[e_1] = [e_2]$, as required.
\end{proof}

By the last proposition, we have $V(C(\mathbb{T},A)) \iso K_0(C(\mathbb{T},A))_+$.
Using the construction of $K_*(A)$ using partial unitaries (see \ccite{8.2.10}{RLL:book}), it is well-known that, for general $C^*$-algebras (not just the ones in Theorem \ref{MainThm}), $K_0(C(\mathbb{T},A))$ is naturally isomorphic to $K_0(A) \dsum K_1(A)$.

Simple, unital, $\jsZ$-stable ASH algebras have stable rank one by \cite{Rordam:Z}.
By \cite{Rieffel:StableRank} and \cite{Brown:StableIsomorphism}, it follows that for every non-zero hereditary subalgebra $B$ of $A \tens \K$, the map from the unitaries of $B^\sim$ to $K_1(A)$ is surjective.
In particular, for any projection $p \in A \tens \K$ and any $[u] \in K_1(A)$, there exists a unitary $v \in p(A \tens \K)p$ such that $[u] = [v+(1-p)]$.
It follows that, for such algebras, $K_0(C(\mathbb{T},A))_+$ can be identified with
\[ \{0\} \cup K_0(A)_{++} \times K_1(A). \]

Next, we show that when $X = \mathbb{T}$, the only non-trivial projection information in Theorem \ref{MainThm} (ie.\ the second coordinate of the invariant) occurs when $f$ takes a constant non-zero value in $V(A)$.

\begin{prop}\alabel{CuT-Descr}
Let $A$ be a $C^*$-algebra.
If $Y \subsetneq \mathbb{T}$ is the intersection of an open and a closed subset of $\mathbb{T}$ and $\langle q \rangle \in V_c(Y,A)$ such that $[q(t)] = [p]$ for all $t \in Y$ then $\langle q \rangle = \langle 1_{C_b(Y)} \tens p \rangle$.
\end{prop}

\begin{proof}
By Lemma \ref{ProjNbhdExtension}, let $U$ be an open neighbourhood of $Y$ upon which $q$ extends continuously, and such that $[q(t)] = [p]$ for all $t \in U$.
Since $Y \neq \mathbb{T}$, we may assume that $U \neq \mathbb{T}$.
Therefore, $U$ decomposes as a (possibly infinite) disjoint union of intervals (given by $\{e^{it}: t \in (a,b)\}$ for some $a < b$ such that $b-a < 2\pi$).
Each of these intervals is contractible, so we can see that the restriction of $q$ to each interval is equivalent to the trivial projection (ie.\ one that is constantly $p$).
It follows that $q$ itself is equivalent to the trivial projection.
\end{proof}
\section{Description of $W(C_0(X,A))$\alabel{OldCu}}

For algebras $A$ as in Theorem \ref{MainThm}, using a fact about the regularity of the Cuntz semigroup for $C_0(X,A)$, we are able to describe the ``classical'' or ``non-stabilized'' Cuntz semigroup,
\[ W(A) := \{[a] \in \Cu(A): a \in M_n(A) \text{ for some } n\}. \]

\begin{cor}\alabel{ClassicalDescription}
Let $A$ be a non-type I, simple, unital ASH algebra with slow dimension growth, and let $X$ be a second countable, locally compact Hausdorff space.
Then $W(C_0(X,A))$ may be identified with pairs $\left(f,\left(\langle q_{[p]} \rangle\right)_{[p] \in V(A)}\right)$, where 
\begin{itemize}
\item $f:X \to W(A)$ is a bounded function which is lower semicontinuous with respect to $\ll$, and 
\item For each $[p] \in V(A)$, $\langle q_{[p]} \rangle$ is an element of $V_c(f^{-1}([p]),A)$ such that $[q_{[p]}(x)] = [p]$ in $V(A)$ for each $x \in f^{-1}([p])$
\end{itemize}
The ordering is given by $\left(f,\left(\langle q_{[p]} \rangle\right)_{[p] \in V(A)}\right) \leq \left(g, \left(\langle r_{[p]} \rangle\right)_{[p] \in V(A)}\right)$ if $f(x) \leq g(x)$ for each $x$, and for each $[p] \in V(A)$,
\[ \langle q_{[p]}|_{f^{-1}([p]) \cap g^{-1}([p])} \rangle = \langle r_{[p]}|_{f^{-1}([p])} \rangle. \]
The addition is given by 
\[ \left(f,\left(\langle q_{[p]} \rangle\right)_{[p] \in V(A)}\right) + \left(g, \left(\langle r_{[p]} \rangle\right)_{[p] \in V(A)}\right) = \left(f+g, \left(\langle s_{[p]} \rangle\right)_{[p] \in V(A)}\right), \]
where for every pair of projections $0 \leq p' \leq p \in A \tens \K$, we have
\[ s_{[p]}|_{f^{-1}([p']) \cap g^{-1}([p-p'])} = q_{[p']} + r_{[p-p']}. \]
(We have that $(f+g)^{-1}([p])$ breaks into disjoint components $f^{-1}([p']) \cap g^{-1}([p-p'])$, and so this definition of $s_{[p]}$ is continuous.)
\end{cor}

\begin{proof}
By \ccite{Theorem 4.5}{Rordam:Z}, we have that the radius of comparison of $\Cu(C_0(X,A))$ is $0$.
Consequently, by \ccite{Theorem 4.12}{BRTTW}, it follows that $W(C_0(X,A))$ is a hereditary subset of $\Cu(C_0(X,A))$.
The description then follows from Theorem \ref{MainThm}, along with the observation that if $f:X \to \Cu(A)$ is strictly bounded by $[a] \in W(A)$ then any element corresponding to a pair $(f,([q_{[p]}]))$ will be bounded by $[1_{C_0(X)} \tens a] \in W(C_0(X,A))$.
\end{proof}

\begin{rmk}
It could be objected that, when we say ``$f$ is lower semicontinuous with respect to $\ll$'', it is not clear whether we mean $\ll$ in $W(A)$ of $\ll$ in $\Cu(A)$ (if $S$ is a subset of an ordered set $T$, $\ll$ may be a different relation when taken with respect to $S$ or with respect to $T$).
However, using \ccite{Theorem 4.12}{BRTTW} and \ccite{Lemma 2.3}{AntoineBosaPerera}, one can see that $\ll$ is the same relation when taken with respect to $W(A)$ or $\Cu(A)$.
\end{rmk}
\section{The Murray-von Neumann semigroup of $C(K,A)$ \alabel{MvN}}
This section contains some remarks on the Murray-von Neumann semigroups that appear in the description of $\Cu(C_0(X,A))$ in Theorem \ref{MainThm}.
The Murray-von \linebreak Neumann semigroups that occur are $V(C(K,A))$ where $K$ is a compact subset of $X$.
Already we have seen, by Proposition \ref{ProjCancellation}, that we have $V(C(K,A)) = K_0(C(K,A))_+$.

The following proposition demonstrates a weak unperforation-type property for $K_0(C(K,A))_+$ (applicable because the algebras $A$ in Theorem \ref{MainThm} are $\jsZ$-stable).

\begin{prop}(cf.\ \ccite{Theorem 1.4}{GongJiangSu:obstructions})\alabel{GJS-Generalized}
Let $A$ be a stably finite, unital $C^*$-algebra and $\iota:A \to A \tens \jsZ$ the canonical embedding.
Let $g \in K_0(A)$.
Then $\iota_*(g) > 0$ if and only if there exists $n_0$ such that $ng > 0$ for all $n \geq n_0$.
\end{prop}

\begin{proof}
This proof is almost identical to that of \ccite{Theorem 1.4}{GongJiangSu:obstructions}.
If $\iota_*(g) > 0$ then, by the standard construction of $\jsZ$, we must have
\[ \iota_{p,q}(g) > 0 \]
in $K_0(A \tens Z_{p,q})$ for some coprime integers $p,q$, where $\iota_{p,q}:A \to A \tens Z_{p,q}$ is the canonical embedding.
In particular, if $\pi_0:A \tens Z_{p,q} \to A \tens M_p$ and $\pi_1:A \tens Z_{p,q} \to A \tens M_q$ denote the evaluation maps at $0$ and $1$ respectively, then we have
\[ qg = (\pi_0 \circ \iota_{p,q})_*(g) > 0 \]
and
\[ pg = (\pi_1 \circ \iota_{p,q})_*(g) > 0. \]
Since $p$ and $q$ are coprime, it follows that for all $n$ sufficiently large, $ng > 0$.

On the other hand, suppose that for all $n \geq n_0$, we have $ng > 0$.
Then, letting $p,q \geq n_0$ be coprime integers, we have for some $e,f \in A \tens \K$ that
\[ pg = [e] \text{ and } qg = [f], \]
and obviously $q[e] = p[f]$.
By \ccite{Theorem 3.1.4}{Blackadar:rational}, we must have $e \tens 1_{qk} \mvneq f \tens 1_{pk}$ for all $k \geq k_0$, for some $k_0$.
Let $k \geq k_0$ be such that $k$ and $q$ are coprime.

Since $e \tens 1_{qk} \mvneq f \tens 1_{pk}$, there exists a homotopy of projections, $E \in C([0,1],A \tens \K)$ such that 
\[ E(0) = e \tens 1_{qk} \text{ and } E(1) = f \tens 1_{pk}. \]
Viewing $E$ as a projection in $A \tens Z_{pk,q} \tens \K$, we can see by \ccite{Lemma 1.2}{GongJiangSu:obstructions} that $[E] = (\iota_{pk,q})_*(g)$.
Hence, $(\iota_{pk,q})_*(g) > 0$ and, since $\iota$ factors through $\iota_{pk,q}$, we have $\iota_*(g) > 0$.
\end{proof}

\begin{rmk}\alabel{Z-MvNdesc}
Combining Propositions \ref{ProjCancellation} and \ref{GJS-Generalized}, we find that $V(C(K,\jsZ))$ can be identified with
\[ \{g \in K_0(C(K)):\ \exists n_0 \text{ s.t.\ } ng \geq 0\ \forall n \geq n_0\}. \]
\end{rmk}

We can also produce a fairly explicit description of $V(C(K,A))$ when $A$ is an AF algebra, by the following result.

\begin{prop}\alabel{AF-K0}
Let $A, B$ be $C^*$-algebras such that $A$ is AF and $B$ is unital.
Then $K_0(A \tens B)$ can be identified with $K_0(A) \tens_{\mathbb{Z}} K_0(B)$, and such that $K_0(A \tens B)_+$ is identified with the ordered cone consisting of sums of elementary tensors $g \tens h$ where $g \in K_0(A)_+$ and $h \in K_0(B)_+$.
\end{prop}

\begin{proof}
It is easy to see that the statement holds in the special case that $A$ is a finite-dimensional $C^*$-algebra; that is, if $A$ is finite-dimensional with $k$ simple direct summands then
\[ K_0(A \tens B) \iso \mathbb{Z}^k \tens K_0(B). \]

In the general case, let us write $A$ as an inductive limit
\[ A_1 \labelledrightarrow{\phi_1^2} A_2 \labelledrightarrow{\phi_2^3} \cdots \rightarrow A, \]
where each algebra $A_i$ is finite-dimensional.
Consequently, we have a direct limit of ordered abelian groups
\[ K_0(A_1 \tens B) \labelledrightarrow{(\phi_1 \tens \id_B)_*} K_0(A_2 \tens B) \rightarrow \cdots \rightarrow K_0(A \tens B), \]
which can be rewritten
\[ K_0(A_1) \tens K_0(B) \labelledrightarrow{(\phi_1)_* \tens \id_{K_0(B)}} K_0(A_2) \tens K_0(B) \rightarrow \cdots \rightarrow K_0(A \tens B). \]

It follows by \ccite{Lemma 2.2}{GoodearlHandelman:tensor} that, as ordered groups (using the tensor product ordering given in the statement of the proposition), we have
\[ K_0(A \tens B) = (\varinjlim K_0(A_i)) \tens K_0(B) = K_0(A) \tens K_0(B). \]
\end{proof}

\begin{rmk}\alabel{AF-MvNdesc}
Combining Propositions \ref{ProjCancellation} and \ref{AF-K0}, we find that when $A$ is AF, we can express
\[ V(C(K,A)) \iso K_0(C(K))_+ \tens K_0(A)_+, \]
as a tensor product of semigroups (as defined in \cite{Grillet}).
Explicitly, this is the abelian semigroup generated by simple tensors $g \tens h$ where $g \in K_0(C(K))_+$ and $h \in K_0(A)_+$, with the usual identifications
\[ (g_1 + g_2) \tens h = g_1 \tens h + g_2 \tens h \]
and
\[ g \tens (h_1 + h_2) = g \tens h_1 + g \tens h_2. \]
\end{rmk}


\begin{thebibliography}{10}

\bibitem{AntoineBosaPerera}
R.~Antoine, J.~Bosa, and F.~Perera.
\newblock Completions of monoids with applications to the {C}untz semigroup.
\newblock {\em Internat. J. Math.},
\newblock to appear. arXiv preprint math.OA/1003.2874.

\bibitem{Blackadar:rational}
B.~Blackadar.
\newblock Rational {$C^*$}-algebras and nonstable {$K$}-theory.
\newblock {\em Rocky Mountain J. Math.}, \textbf{20}(2) (1990), 285--316.

\bibitem{Blackadar:KBook}
B.~Blackadar.
\newblock {\em {$K$}-theory for operator algebras}, volume~5 of {\em
  Mathematical Sciences Research Institute Publications}.
\newblock Cambridge University Press, Cambridge, second edition, 1998.

\bibitem{BlackadarHandelman}
B.~Blackadar and D.~Handelman.
\newblock Dimension functions and traces on {$C^{\ast} $}-algebras.
\newblock {\em J. Funct. Anal.}, \textbf{45}(3) (1982), 297--340.

\bibitem{BRTTW}
B.~Blackadar, L.~Robert, A.~P. Tikuisis, A.~S. Toms, and W.~Winter.
\newblock An algebraic approach to the radius of comparison.
\newblock arXiv preprint. math.OA/1008.4024, Aug. 2010.

\bibitem{Brown:StableIsomorphism}
L.~G. Brown.
\newblock Stable isomorphism of hereditary subalgebras of {$C\sp*$}-algebras.
\newblock {\em Pacific J. Math.}, \textbf{71}(2) (1977), 335--348.

\bibitem{BrownCiuperca}
N.~P. Brown and A.~Ciuperca.
\newblock Isomorphism of {H}ilbert modules over stably finite {$C\sp
  *$}-algebras.
\newblock {\em J. Funct. Anal.}, \textbf{257}(1) (2009), 332--339.

\bibitem{BrownPereraToms}
N.~P. Brown, F.~Perera, and A.~S. Toms.
\newblock The {C}untz semigroup, the {E}lliott conjecture, and dimension
  functions on {$C\sp *$}-algebras.
\newblock {\em J. Reine Angew. Math.}, \textbf{621} (2008), 191--211.

\bibitem{BrownWinter:quasitraces}
N.~P. Brown and W.~Winter.
\newblock Quasitraces are traces: A short proof of the finite-nuclear-dimension
  case.
\newblock arXiv preprint. math.OA/1005.2229, May 2010.

\bibitem{CiupercaElliottSantiago}
A.~Ciuperca, G.~Elliott, and L.~Santiago.
\newblock On inductive limits of type {I} {$C^*$}-algebras with one-dimensional
  spectrum.
\newblock {\em Int. Math. Res. Not. IMRN},
\newblock to appear. arXiv preprint math.OA/1004.0262.

\bibitem{CiupercaRobertSantiago}
A.~Ciuperca, L.~Robert, and L.~Santiago.
\newblock The {C}untz semigroup of ideals and quotients and a generalized
  {K}asparov stabilization theorem.
\newblock {\em J. Operator Theory}, \textbf{64}(1) (2010), 155--169.

\bibitem{CowardElliottIvanescu}
K.~T. Coward, G.~A. Elliott, and C.~Ivanescu.
\newblock The {C}untz semigroup as an invariant for {$C\sp *$}-algebras.
\newblock {\em J. Reine Angew. Math.}, \textbf{623} (2008), 161--193.

\bibitem{Cuntz:dimensionfunctions}
J.~Cuntz.
\newblock Dimension functions on simple {$C\sp*$}-algebras.
\newblock {\em Math. Ann.}, \textbf{233}(2) (1978), 145--153.

\bibitem{Edwards}
D.~A. Edwards.
\newblock S\'eparation des fonctions r\'eelles d\'efinies sur un simplexe de
  {C}hoquet.
\newblock {\em C. R. Acad. Sci. Paris}, \textbf{261} (1965), 2798--2800.

\bibitem{ElliottRobertSantiago}
G.~Elliott, L.~Robert, and L.~Santiago.
\newblock The cone of lower semicontinuous traces on a {$C^*$}-algebra.
\newblock {\em Amer. J. Math.}, 
\newblock to appear. arXiv preprint math.OA/0805.3122.

\bibitem{ElliottGongLi:AHclassification}
G.~A. Elliott, G.~Gong, and L.~Li.
\newblock On the classification of simple inductive limit {$C^*$}-algebras.
  {II}. {T}he isomorphism theorem.
\newblock {\em Invent. Math.}, \textbf{168}(2) (2007), 249--320.

\bibitem{GongJiangSu:obstructions}
G.~Gong, X.~Jiang, and H.~Su.
\newblock Obstructions to {$\mathcal Z$}-stability for unital simple
  {$C^*$}-algebras.
\newblock {\em Canad. Math. Bull.}, \textbf{43}(4) (2000), 418--426.

\bibitem{Goodearl:book}
K.~R. Goodearl.
\newblock {\em Partially ordered abelian groups with interpolation}, volume~20
  of {\em Mathematical Surveys and Monographs}.
\newblock American Mathematical Society, Providence, RI, 1986.

\bibitem{GoodearlHandelman:tensor}
K.~R. Goodearl and D.~E. Handelman.
\newblock Tensor products of dimension groups and {$K_0$} of unit-regular
  rings.
\newblock {\em Canad. J. Math.}, \textbf{38}(3) (1986), 633--658.

\bibitem{Grillet}
P.-A. Grillet.
\newblock The tensor product of semigroups.
\newblock {\em Trans. Amer. Math. Soc.}, \textbf{138} (1969), 267--280.

\bibitem{Haagerup:quasitraces}
U.~Haagerup.
\newblock Quasi-traces on exact {$C^*$}-algebras are traces.
\newblock Unpublished preprint, 1991.

\bibitem{LinNiu:Lifting}
H.~Lin and Z.~Niu.
\newblock Lifting {$KK$}-elements, asymptotic unitary equivalence and
  classification of simple {$C^\ast$}-algebras.
\newblock {\em Adv. Math.}, \textbf{219}(5) (2008), 1729--1769.

\bibitem{Lin:AsympClassification}
H.~X. Lin.
\newblock Asymptotically unitary equivalence and classification of simple
  amenable {$C^*$}-algebras.
\newblock {\em Invent. Math.}, 
\newblock to appear. arXiv preprint math.OA/0806.0636.

\bibitem{LinPhillips:LimitDecomp}
Q.~Lin and N.~C. Phillips.
\newblock Direct limit decomposition for {$C^*$}-algebras of minimal
  diffeomorphisms.
\newblock In {\em Operator algebras and applications}, volume~38 of {\em Adv.
  Stud. Pure Math.}, pages 107--133. Math. Soc. Japan, Tokyo, 2004.

\bibitem{NgWinter:ash}
P.~W. Ng and W.~Winter.
\newblock A note on subhomogeneous {$C^\ast$}-algebras.
\newblock {\em C. R. Math. Acad. Sci. Soc. R. Can.}, \textbf{28}(3) (2006),
  91--96.

\bibitem{PereraToms:recasting}
F.~Perera and A.~S. Toms.
\newblock Recasting the {E}lliott conjecture.
\newblock {\em Math. Ann.}, \textbf{338}(3) (2007), 669--702.

\bibitem{Phillips:Tori}
N.~C. Phillips.
\newblock Every simple higher dimensional noncommutative torus is an
  {A$\mathbb{T}$} algebra.
\newblock arXiv preprint. math.OA/0609783, Sept. 2006.

\bibitem{Phillips:arsh}
N.~C. Phillips.
\newblock Cancellation and stable rank for direct limits of recursive
  subhomogeneous algebras.
\newblock {\em Trans. Amer. Math. Soc.}, \textbf{359}(10) (2007), 4625--4652
  (electronic).

\bibitem{Phillips:rsh}
N.~C. Phillips.
\newblock Recursive subhomogeneous algebras.
\newblock {\em Trans. Amer. Math. Soc.}, \textbf{359}(10) (2007), 4595--4623
  (electronic).

\bibitem{Rieffel:StableRank}
M.~A. Rieffel.
\newblock Dimension and stable rank in the {$K$}-theory of
  {$C^{\ast}$}-algebras.
\newblock {\em Proc. London Math. Soc. (3)}, \textbf{46}(2) (1983), 301--333.

\bibitem{Robert:NCCW}
L.~Robert.
\newblock Classification of inductive limits of 1-dimensional {NCCW} complexes.
\newblock arXiv preprint. math.OA/1007.1964, July 2010.

\bibitem{CommutativeCuntz}
L.~Robert and A.~Tikuisis.
\newblock Hilbert {$C^*$}-modules over a commutative {$C^*$}-algebra.
\newblock {\em Proc. Lond. Math. Soc. (3)}, 
\newblock to appear. arXiv preprint math.OA/0910.2967.

\bibitem{Rordam:UHFII}
M.~R{\o}rdam.
\newblock On the structure of simple {$C^*$}-algebras tensored with a
  {UHF}-algebra. {II}.
\newblock {\em J. Funct. Anal.}, \textbf{107}(2) (1992), 255--269.

\bibitem{Rordam:Z}
M.~R{\o}rdam.
\newblock The stable and the real rank of {$\mathcal Z$}-absorbing
  {$C^*$}-algebras.
\newblock {\em Internat. J. Math.}, \textbf{15}(10) (2004), 1065--1084.

\bibitem{RLL:book}
M.~R{\o}rdam, F.~Larsen, and N.~Laustsen.
\newblock {\em An introduction to {$K$}-theory for {$C^*$}-algebras}, volume~49
  of {\em London Mathematical Society Student Texts}.
\newblock Cambridge University Press, Cambridge, 2000.

\bibitem{Thomsen:Solenoids}
K.~Thomsen.
\newblock The homoclinic and heteroclinic {$C^*$}-algebras of a generalized
  one-dimensional solenoid.
\newblock {\em Ergodic Theory Dynam. Systems}, \textbf{30}(1) (2010), 263--308.

\bibitem{Toms:rigidity}
A.~S. Toms.
\newblock {K}-theoretic rigidity and slow dimension growth.
\newblock {\em Invent. Math.}, 
\newblock to appear. arXiv preprint math.OA/0910.2061.

\bibitem{Toms:comparison}
A.~S. Toms.
\newblock Comparison theory and smooth minimal {$C\sp *$}-dynamics.
\newblock {\em Comm. Math. Phys.}, \textbf{289}(2) (2009), 401--433.

\bibitem{TomsWinter:ssa}
A.~S. Toms and W.~Winter.
\newblock Strongly self-absorbing {$C^*$}-algebras.
\newblock {\em Trans. Amer. Math. Soc.}, \textbf{359}(8) (2007), 3999--4029.

\bibitem{TomsWinter:MinClassification}
A.~S. Toms and W.~Winter.
\newblock Minimal dynamics and the classification of {$C^*$}-algebras.
\newblock {\em Proc. Natl. Acad. Sci. USA}, \textbf{106}(40) (2009),
  16942--16943.

\bibitem{Winter:drZstable}
W.~Winter.
\newblock Decomposition rank and {$\mathcal Z$}-stability.
\newblock {\em Invent. Math.}, \textbf{179}(2) (2010), 229--301.

\bibitem{Winter:Perfect}
W.~Winter.
\newblock Nuclear dimension and $\mathcal{Z}$-stability of perfect
  {C}$^*$-algebras.
\newblock arXiv preprint. math.OA/1006.2731, June 2010.

\end{thebibliography}
\end{document}